\documentclass[12pt]{amsart}

\usepackage{amsthm}

\usepackage{graphicx}
\usepackage{epsfig}
\usepackage{epstopdf}

\usepackage{latexsym}
\usepackage{amsmath}
\usepackage{amssymb}
\usepackage{amsfonts}

\theoremstyle{definition}

\newtheorem{dfn}{Definition}[section]
\newtheorem{defn}[dfn]{Definition}
\newtheorem{definition}[dfn]{Definition}

\newtheorem{example}[dfn]{Example}

\newtheorem{remark}[dfn]{Remark}
\theoremstyle{plain}
\newtheorem*{ack}{Acknowledgments} 
\newtheorem{problem}[dfn]{Problem}

\newtheorem{theorem}[dfn]{Theorem}
\newtheorem{lem}[dfn]{Lemma}
\newtheorem{lemma}[dfn]{Lemma}

\newtheorem{proposition}[dfn]{Proposition}

\newtheorem{assumption}[dfn]{Assumption}
\newtheorem{cor}[dfn]{Corollary}

\newtheorem{conjecture}[dfn]{Conjecture}

\newtheorem{question}[dfn]{Question}

\newtheorem{notation}[dfn]{Notation}

\def\proof{\par\medskip\noindent{\it Proof. }}

\def\CB{\mathcal B}
\def\CC{\mathcal C}
\def\CD{\mathcal D}
\def\CE{\mathcal E}
\def\CK{\mathcal K}
\def\CL{\mathcal L}
\def\CO{\mathcal O}

\def\CP{\mathcal P}
\def\CV{\mathcal V}

\def\VV{\mathbf V}

\def\P{{\mathbb P}}

\def\R{{\mathbb R}}
\def\C{{\mathbb C}}
\def\Z{{\mathbb Z}}

\def\H{{\mathbb H}}
\def\N{{\mathbb N}}

\def\al{\alpha}
\def\be{\beta}
\def\ga{\gamma}
\def\Ga{\Gamma}

\def\Si{\Sigma}
\def\si{\sigma}

\def\la{\lambda}
\def\La{\Lambda}
\def\Om{\Omega}

\def\ul{\underline}
\def\acts{\curvearrowright}
\def\D{\partial}

\def\embed{\hookrightarrow}
\def\<{\langle}
\def\>{\rangle}
\def\geo{\partial_{\infty}}

\def\rank{\mathop{\hbox{rank}}}

\def\t{\tilde}

\newcommand{\mor}[0]{\operatorname{Mor}} 

\newcommand{\Span}{\operatorname{Span}}
\newcommand{\Faces}{\operatorname{Faces}}
\newcommand{\Resid}{\operatorname{Res}}
\newcommand{\Ob}{\operatorname{Ob}}
\newcommand{\Nerve}{\operatorname{Nerve}}
\newcommand{\inter}[0]{\operatorname{Int}}    
\newcommand{\im}{\operatorname{im}}

\begin{document}

\title{Dirichlet fundamental domains and complex--projective varieties} 
\author{Michael Kapovich}
\date{\today}

\begin{abstract}
We prove that for every finitely-presented group $G$  there exists a $2$-dimensional irreducible complex-projective variety $W$ with the fundamental group $G$, so that all singularities of $W$ are normal crossings and Whitney umbrellas. 
\end{abstract}

\maketitle

\section{Introduction}

It is well-known that fundamental groups of compact K\"ahler manifolds satisfy many restrictions, see e.g. \cite{ABC}. On the other hand, C.~Simpson 
proved in \cite{Simpson} that every finitely-presented group $G$ appears as the fundamental group of a (singular) irreducible complex-projective variety. In the same paper Simpson asked the following question which is a variation on a problem about fundamental groups of irreducible projective varieties originally posed by D.~Toledo: 

\begin{question}\label{ques:simpson}
Is it true that every finitely-presented group $G$ is isomorphic to the fundamental group of a irreducible complex-projective variety 
whose singularities are normal crossings only? 
\end{question}

In our previous paper with J\'anos Koll\'ar 
\cite{KK} we proved that the answer to this question is positive provided one does not require irreducibility. Although we do not 
know what the answer to Simpson's original question is, in this paper we prove 

%Our main result is

\begin{theorem}\label{thm:main}
Let $G$ be a finitely-presented group. Then there exists a $2$-dimensional irreducible complex-projective variety $W$ with the fundamental group $G$, so that the only singularities of $W$ are normal crossings and Whitney umbrellas. 
%$\Z_2$-quotients of normal crossing singularities, so that the local model for the action of $\Z_2=\<\theta\>$ is
%$$y_1  y_2=0, \quad \theta(y_1,y_2,y_3,y_4)=  (y_2,y_1, y_3, -y_4).  $$
Furthermore, if $G$ is isomorphic to the fundamental group of a compact 3-dimensional hyperbolic manifold with (possibly empty) 
convex boundary, then all singularities of $W$ are normal crossings. 
\end{theorem}

In other words, we get $W$ with ``controlled'' singularities (unlike the ones which appear in Simpson's proof in \cite{Simpson}). 
The key tools in our proof are the recent ``universality'' theorem  by Petrunin and Panov, and a certain genericity result for Dirichlet fundamental domains of discrete isometry groups of the hyperbolic 3-space. 

\begin{theorem}
(D.Panov, A.Petrunin, \cite{PP}) 
\label{T1.0}
Let $G$ be a finitely-presented group. Then there exists a discrete cocompact 
subgroup $\Gamma< PO(3,1)$ so that:

1. The only nontrivial finite subgroups of $\Gamma$ are isomorphic to $\Z_2$ or $\Z_2\times \Z_2$. 

2. For each order 2 element of $\Ga$ its fixed-point set in the hyperbolic 3-space has dimension 0 or 1. 

3. The fundamental group of the quotient $M:=\H^3/\Gamma$ is  isomorphic to $G$. 
\end{theorem}

Note that $M$ is a 3-dimensional complex which is a manifold away from a finite subset, where the singularities are cones over projective planes. We will need a minor variation on their construction:

\begin{theorem}\label{T1.1}
Let $G$ be a finitely-presented group. Then there exists a discrete  nonelementary 
subgroup $\t\Gamma< PO(3,1)$ so that: 

1. All nontrivial finite subgroups of $\t\Gamma$ are isomorphic to $\Z_2$, each has a single fixed point in $\H^3$. 
In other words, every nontrivial finite  subgroup of $\t\Gamma$ is generated by a Cartan involution of $\H^3$. 

2. The group $\t\Ga$ is convex-cocompact (every convex fundamental domain in $\H^3$ 
of $\t\Ga$ has only finitely many faces and $\t\Ga$ contains no parabolic elements). 

3. The fundamental group of the quotient $\H^3/\t\Gamma$ is  isomorphic to $G$.
\end{theorem}

We will refer to the class of subgroups of $PO(3,1)$ satisfying property 1 in this theorem as {\em class} $\CK$ and to the class 
of groups satisfying properties 1 and 2 as the {\em class} $\CK^2$.

For a discrete subgroup $\Gamma< PO(3,1)$ and a point $x\in \H^3$ (not fixed by any nontrivial element of $\Gamma$) we define 
the {\em Dirichlet tiling} $\CD_x$ of $\H^3$ to be the Voronoi tiling of $\H^3$ corresponding to the orbit 
$\Gamma \cdot x$. The {\em tiles} of $\CD_x$ are the Dirichlet fundamental domains
$$
D_{\ga x}=\{p\in \H^3: d(p,\ga x)\le d(p, \al(x)), \forall \al \in\Gamma\setminus \{\ga\} \}
$$

\begin{conjecture}\label{conj:C1}
For $\Ga<PO(3,1)$ of class $\CK$, for generic choice of $x$ the tiling ${\mathcal D}_x$ {\em simple}, i.e., the dual cell-complex 
to  ${\mathcal D}_x$  is a simplicial complex. 
\end{conjecture}

Conjecture \ref{conj:C1} was stated as a theorem (for torsion-free groups $\Gamma$) in the paper by Jorgensen and Marden \cite{JM}. 
However, their proof has a serious gap noted by Diaz and Ushijima in 
\cite{DU}: The trouble with \cite{JM} is confusion between algebraic and semi-algebraic sets. In \cite{JM} one of the key claims (Corollary 3.1) is that certain semi-algebraic sets in $\H^3$ have empty interiors,  while all what they proved  is that these are  
proper subsets of $\H^3$. (The sets in question are subsets $\CE({\ul{A}})\subset \H^3$ consisting of points $x$ such that quadruple intersections 
of bisectors 
$$
\bigcap_{i=1}^4 Bis(x, A_i x),
$$
are non-transversal in $\H^3$. Here $\ul{A}=\{A_1,...,A_4\}$, where $A_i\in \Ga$ are fixed pairwise distinct and nontrivial elements.) 
We will actually see in \S \ref{sec:examples}  that some of the sets $\CE({\ul{A}})$ could have non-empty interiors. 
The paper \cite{DU} proves an analogue of Conjecture \ref{conj:C1} for torsion-free orientation-preserving discrete subgroups of $PO(2,1)$. We 
do not know how to prove Conjecture \ref{conj:C1} either. Nevertheless, we will prove a weaker result that will suffice for our purposes:

\begin{theorem}\label{thm:generic} 
Suppose that $\Ga<PO(3,1)$ is a subgroup of class $\CK$. Then for a generic choice of $x\in \H^3$ the Dirichlet 
tiling ${\mathcal D}_x$ is simple away from its vertex set $\CD_x^{(0)}$. 
Moreover, only points in the interiors of 2-dimensional faces of $\CD_x$ 
can be fixed by Cartan involutions in $\Ga$. 
\end{theorem}

\begin{remark}
After completing this paper I received the preprint \cite{U} by Ushijima where Conjecture \ref{conj:C1} is proven for purely loxodromic subgroups 
of $PO(3,1)_+\cong PSL(2,\C)$.  The arguments in \cite{Ushijima} are different from the ones used in the proof of Theorem \ref{thm:generic}. 
\end{remark}

Once Theorems  \ref{T1.1} and \ref{thm:generic}  are established, the proof of Theorem \ref{thm:main} follows closely 
the arguments in \cite{KK}, by complexifying a certain hyperbolic polyhedral complex $\CC$ (obtained by taking a quotient of 
${\mathcal D}_x\setminus \CD_x^{(0)}$) and then blowing up ``parasitic subspaces'' of the complexification. Using $\CC$ one constructs a (reducible) projective variety $X$ and a finite group $\Theta$ acting on $X$, so that the only singularities of $X$ are normal crossings, 
the projective variety $V= X/\Theta$ is irreducible and $\pi_1(V)\cong G$. Irreducibility 
of $V$ comes from the fact that all facets of ${\mathcal D}_x$ are equivalent under the 
$\tilde\Ga$-action  (unlike the Euclidean polyhedral complexes used in \cite{KK} which have many facets). 
The projective surface $W$ is obtained by applying Lefschetz hyperplane 
section theorem to $V$. Whitney umbrella singularities of $W$ correspond to the fixed 
points of the action on $\C \P^3$ of Cartan involutions in $\tilde\Ga$.  

 \begin{ack} {\em This paper grew out of our work \cite{KK} with J\'anos Koll\'ar and 
 I am  grateful to him for questions, comments and suggestions. In particular, he explained to me that 
 Whitney umbrella singularities appear as $\Z_2$-quotients of normal crossings and suggested the 
 dimension reduction from 3 to 2. I am also grateful to Akira Ushijima for sharing with me \cite{Ushijima}.  
 Partial financial support for this work was provided by the NSF 
 grant number DMS-09-05802.}
\end{ack}

\section{Preliminaries}\label{sec:2}

\begin{notation}
Throughout the paper we will use the topologist's convention: \newline $\Z_2=\Z/2\Z$. 
\end{notation}

\medskip
Let $\R^{n,1}$ denote the Lorentzian space, it is $\R^{n+1}$ equipped with nondegenerate inner product $x\cdot y$ of the signature 
$(n,1)$. We will be mostly interested in the case $n=3$, but our proofs are more general. 
We will refer to the inner product $x\cdot x$ as the {\em Lorentzian norm} of $x$. 
The {\em light cone} ${\CL}$ of $\R^{n,1}$ consists of vectors of negative Lorentzian norm. This cone 
has two components, we fix one of these components  ${\mathcal L}^\uparrow$; we will refer to 
${\mathcal L}^\uparrow$ as the {\em future light cone}. We let $C$ denote the boundary of  ${\mathcal L}^\uparrow$ 
and $C^+$ the closure of ${\mathcal L}^\uparrow$. The {\em future} (or the ``upper'') sheet of the hyperboloid 
$$
\{x| x\cdot x=-1\}
$$
is the intersection $H$ of this hyperboloid with ${\mathcal L}^\uparrow$. Then $H$ is the Lorentzian model of the hyperbolic 
$n$-space $\H^n$: Restriction of the Lorentzian inner product to the tangent bundle of $H$ is a Riemannian metric 
of the sectional curvature $-1$ on $H$. For a subset $E\subset \R^{n+1}$ we let $\P E$ denote its projection to $\R \P^n$.  
We will identify $\H^n$ with the projectivization $\P H$ of $H$ (and of ${\mathcal L}^\uparrow$). The projectivization $\P C^+$ 
of the cone $C^+$ is the standard compactification of $\H^n$: $\P C^+= \H^n \cup S^{n-1}$, where $S^{n-1}=\P C$.  
For a subset $X$ of $\H^n$ we define its {\em ideal boundary} $\geo X$ by:
$$
\P(cl(X)\cap C).
$$
In other words $\geo X$ is the accumulation of $X$ on the boundary sphere $S^{n-1}$ of $\H^n$.

For $x,y\in H$ we let $d(x,y)$ denote their hyperbolic distance. Then (see e.g. \cite{Ratcliffe}) 
\begin{equation}\label{cdot}
x\cdot y= -\cosh(d(x,y)). 
\end{equation}
In particular, $x\cdot y=-1$ iff $x=y$.

\begin{lemma}\label{L2.0}
Let $u, v, w\in H$. 
%be vectors of the same nonzero Lorentzian norm. 
Then for any $t,s\in \R$ such that $st\ne 0, s+t\ne 0$,   
$$
su+ tv \ne (s+t)w 
$$
unless $u=v=w$. 
\end{lemma}
\proof Note that it suffices to show that $u=v$ (since $s+t\ne 0$). 
%We assume that $u, v, w$, are such that 
%$$u\cdot u=v\cdot v=w\cdot w=-1. $$
%the case when their inner products equal to $1$ are similar. 
%Suppose that $tu+ sv = (t+s)w$. 
Computing the Lorentzian norms of both sides of the equation 
$$su+ tv = (s+t)w$$
 we get:
$$
uv=-1. 
$$
Since $u, v\in H$, 
%Since for $u, v\in H$, $u\cdot v= -\cosh(d(u,v))$, where $d$ is the hyperbolic distance  between $u, v$, 
it follows that $u=v$. Since $t+s\ne 0$, it follows that $u=v=w$. \qed 
%Note that we cannot have the situation when $u\in C^+$ and $v\notin C^+$ for then $u\cdot v\ge 1$. \qed 

\medskip
For $x, y\in \H^n$ the {\em bisector} $Bis(x,y)$ is the hyperplane
$$
Bis(x,y)=\{p\in \H^n: d(x,p)=d(y,p)\}
$$
In view of the equation \eqref{cdot}, bisectors are described by
$$
Bis(x,y)=\{p\in H: x\cdot p=y\cdot p\}.  
$$
We extend this definition to the entire $\R^{n,1}$, then the {\em extended} bisector $\widetilde{Bis}(x,y)$ is the hyperplane
$$
\widetilde{Bis}(x,y)=\{p\in \R^{n+1}: p\cdot (x-y)=0\}= (x-y)^\perp. 
$$
The (extended) bisectors $\widetilde{Bis}(x_i, y_i), i=1,...,k$, are transversal in $\R^{n+1}$ iff the normal vectors $(x_i-y_i), i=1,...,k$ 
are linearly independent. In particular, in order to verify transversality of the bisectors 
${Bis}(x_i, y_i), i=1,...,k$ in $\H^n$, it suffices to check linear independence of the vectors 
$(x_i-y_i), i=1,...,k$.

\medskip
{\bf Isometry group.} We let $O(n,1)<GL(n,\R)$ denote the automorphism group of $\R^{n,1}$. This group has index 2 subgroup 
$O(n,1)^{\uparrow}$ preserving the future light cone ${\mathcal L}^\uparrow$. Thus, 
$$
O(n,1)=O(n,1)^{\uparrow} \times \Z_2
$$
where $\Z_2=\{\pm I\}$ and $I\in GL(n+1,\R)$ is the identity matrix. In particular, 
$O(n,1)^{\uparrow}$ is isomorphic to $PO(n,1)$, the isometry group of $\H^n$. This isomorphism will allow us to 
identify subgroups of $PO(n,1)$ with subgroups of $O(n,1)$. 

\medskip
{\bf Classification of nontrivial elements of $PO(n,1)$.} 

1. An element $A\in O(n,1)^{\uparrow}$ is {\em elliptic} if it has 
a fixed vector in $H$. If this fixed vector is unique and $n=3$, then $A$ necessarily has order $2$ and reverses orientation. 
Regarding $H=\H^3$ as a symmetric space, such elliptic elements are {\em Cartan involutions} in $\H^3$. For arbitrary $n$, Cartan involutions are characterized by the property that each has a unique fixed point in $\H^n$ (and order $2$). If $A\in PO(n,1)$ has finite order, it is necessarily elliptic. 

\medskip 
2. An element $A\in O(n,1)^{\uparrow}$ is {\em parabolic} if it has a unique, up to a multiple, (nonzero) fixed vector $p\in C^+$ 
and, furthermore, $p$ belongs to $C$.  

\medskip
3. The rest of the isometries of $\H^n$ are {\em loxodromic}. These elements $A\in O(n,1)^{\uparrow}$ are characterized by the property 
that each has exactly two (up to multiple) eigenvectors $e_+, e_-$ in $C$ and the corresponding eigenvalues $\la, \la^{-1}$ are different from $1$. The span of these eigenvectors is a plane $E_{A}\subset \R^{n+1}$ invariant under $A$. The intersection $E_{A}\cap H$ 
is a hyperbolic geodesic $L$ invariant under $A$, it is called the {\em axis} of $A$. The restriction of $A$ to its axis is a nontrivial translation. The eigenvectors $e_\pm$ project to the points in $\P C\cong S^{n-1}$ fixed by $A$. 
(These are the only fixed points that $A$ can have.)

There is a finer classification of loxodromic isometries. Every loxodromic $A\in O(n,1)^{\uparrow}$  preserves 
the orthogonal complement $E_{A}^\perp$ of $E_{A}$. The restriction of the Lorentzian inner product to $E_{A}^\perp$ is 
necessarily positive-definite. If $A$ fixes $E_{A}^\perp$ pointwise, it is called {\em hyperbolic}; otherwise, it is called 
{\em strictly loxodromic}. For $n=3$, this can be described more precisely: 
Each strictly loxodromic element acts on $E_{A}^\perp$ as a nontrivial rotation $R_\theta$ 
(by the angle $\theta$) or a reflection (in case $A$ reverses orientation on $H$). The angle $\theta$ is the {\em angle of rotation} of 
$A$. Intrinsically, in terms of the geometry of $\H^n$, hyperbolic isometries are characterized as compositions $\tau_1 \circ \tau_2$ 
of distinct Cartan involutions. It we trivialize the normal bundle of the axis $L$ by parallel vector fields, then each 
hyperbolic isometry acts on the normal bundle as the translation along the axis, while a strictly loxodromic element is a composition 
of the translation and a (nontrivial) orthogonal transformation of a normal plane. 

\medskip
{\bf Discrete subgroups.} Suppose that $\Ga<PO(n,1)$ is a discrete subgroup. If $A_1, A_2\in \Ga$ are loxodromic which share a common fixed point in $S^{n-1}= \P C$, i.e., they have a common eigenvector $e_+\in C$. 
Then $A_1, A_2$ share the other eigenvector $e_-$ in $C$ as well,  
and, hence,  have the common axis $L$ in $\H^n$, see e.g. \cite{Ratcliffe}. 
If $A_1, A_2\in SO(3,1)$ then they necessarily commute in this situation, 
otherwise, they (typically) do not commute. However, if they do not commute, then the group $\<A_1, A_2\>$ generated by $A_1, A_2$ 
does not act faithfully on $L$, i.e., it contains an elliptic element fixing $L$ pointwise. Thus, such pairs of loxodromic elements cannot belong to a group $\Ga< PO(n,1)$ of the class $\CK$.     

\medskip 
{\bf Fundamental domains.} Let $D$ be a closed 
convex domain in $\H^n$ and $\Ga$ a discrete group of isometries of $\H^n$. Then $D$ is said to be a {\em fundamental domain} 
of $\Ga$ if the following hold:

\begin{enumerate}
\item $\Ga\cdot D=\H^n$. 

\item For every $\ga\in \Ga\setminus \{1\}$, $\ga D\cap D\ne \emptyset$ unless $\ga D\cap D$ is contained in the boundary of $D$. 

\item The covering $\{ \ga D, \ga\in \Ga\}$ of $\H^n$ is {\em locally finite}, i.e., every compact in $\H^n$ intersects only finitely many 
domains $\ga D$. 
\end{enumerate}

In particular, if $D_1, D_2$ are fundamental domains of $\Ga$ and $D_1\subset D_2$ then $D_1=D_2$.

Examples of fundamental domains are provided by the {\em Dirichlet fundamental domain} with the center at $x\in \H^n$, where $x$ 
is not fixed by any $\ga\in \Ga\setminus \{1\}$:
$$
D_x:=\{p\in \H^n: d(p,x)\le d(p, \gamma(x)), \forall \gamma\in \Gamma\}. 
$$
The fundamental domain $D_x$ and its images $D_{\ga x}$ under $\ga\in \Ga$ form a $\Ga$-invariant tiling $\CD_x(\Ga)$ 
of $\H^n$, called the {\em Dirichlet tiling}.% or {\em Dirichlet complex}, since it is a hyperbolic polyhedral complex.  

\medskip
{\bf Convex-cocompact subgroups.} 

A discrete subgroup $\Ga<PO(n,1)$ is called {\em convex-cocompact} if the following holds:

a) $\Ga$ contains no parabolic elements. 

b) %If $D$ is a convex fundamental domain of $\Ga$, then $D$ has only finitely many faces. 
One (equivalently, every) Dirichlet fundamental domain $D_x$ of $\Ga$ is a polyhedral domain in $\H^n$ with finitely many faces.

\medskip
The reader can find detailed discussion and alternative characterizations of convex-cocompact subgroups of $PO(n,1)$ in 
\cite{Beardon-Maskit, Bowditch}. 

{\bf Classes $\CK$ and $\CK^2$.} We say that a subgroup $\Ga<PO(n,1)$ belongs to the class $\CK$ if it is discrete and 
every elliptic element of $\Ga$ is a Cartan involution. We define the class $\CK^2$ to consists of all convex--cocompact discrete subgroups $\Ga< PO(n,1)$ which belong to the class $\CK$.

\medskip
{\bf Elementary and nonelementary groups.} A subgroup $\Ga< PO(n,1)$ is called {\em elementary} if it either has a fixed point 
in the compactification $\H^n\cup S^{n-1}$ of $\H^n$ (this compactification is the projectivization $\P C^+$ of $C^+$) or has an invariant geodesic 
in $\H^n$. Clearly, elementary groups can have nontrivial center, e.g., we can take $\Ga$ to be abelian. Furthermore, 
if $\Ga$ is discrete and elementary then it is {\em virtually abelian}, i.e., contains an abelian subgroup of finite index. 

\begin{lemma}\label{lem:nocenter}
If $\Ga$ is nonelementary and belongs to the class $\CK$, then $\Ga$ has trivial center. 
\end{lemma}
\proof Let $\ga\in \Ga$ be a nontrivial 
central element. Then $\Ga$ preserves the fixed-point set $Fix(\ga)$ of $\ga$ in $\P C^+$. If this fixed set is finite, then 
$\Ga$ is elementary. Otherwise,  $Fix(\ga)$ has to contain a hyperbolic geodesic, so $\ga$ is an elliptic element that cannot be a Cartan involution. 
Contradiction. \qed

\section{Hyperbolic polyhedral complexes}\label{sec:complexes}

This section essentially repeats the definitions given in \cite{KK} in the context of Euclidean polyhedral complexes  
except that we weaken one of the axioms of the polyhedral complex. 

A (convex) {\em hyperbolic polyhedron} is a subset $P$ of $\H^n$ given as the intersection of finitely many open and closed hyperbolic 
half-spaces. Equivalently, if we identify $\H^n$ with the upper sheet $H$ of the hyperboloid in $\R^{n,1}$, then $P$ is given a the 
intersection of $H$ with some convex polyhedral cone $\widetilde{P}\subset \R^{n+1}$: $\widetilde{P}$ is given by a finite set of strict and non-strict linear homogeneous inequalities. Projectivizing $\widetilde{P}$ we obtain a convex projective polytope 
$\widehat{P}\subset \R\P^n$. Note that the correspondence 
$$
P\to \widetilde{P} \to \widehat{P}
$$
is a 1-1 correspondence between convex hyperbolic polyhedra $P\subset \H^n$,  
the convex polyhedral cones $\widetilde{P}\subset \R^{n+1}$ and convex projective polytopes $\widehat{P}\subset \R\P^n$. 
The {\em projective span} $\Span(P)$ is the smallest projective subspace in $\R\P^n$ containing $\widehat{P}$. 
The {\em dimension} of $P$ is its topological dimension, which is the same as the dimension of its projective span $\Span(P)$. 
 
A {\em face} of $P$ is a subset of $P$ which is given by converting some of 
these non-strict inequalities to equalities. Define the set $\Faces(P)$ to be the set of faces of $P$. 
The {\em interior} $\inter(P)$ of $P$ is the topological interior of $P$ in $\Span(P)$. Again, 
$\inter(P)$ is a  hyperbolic polyhedron. We will refer to $\inter(P)$ as an {\em open polyhedron}. 

An (isometric)  {\em morphism} of two hyperbolic polyhedra is an isometric map $f: P\to Q$ so that 
$f(P)$ is  a face of $Q$.  

\begin{definition}
A {\em hyperbolic polyhedral complex}  is a small category $\CC$ whose objects are convex hyperbolic polyhedra and 
morphisms are their isometric morphisms satisfying the following axioms:

Axiom 1. For every $c_1\in \Ob(\CC)$ and every face $c_2$ of $c_1$, $c_2\in \Ob(\CC)$,  the inclusion map $\iota: c_1\to c_2$ is a morphism of 
$\CC$.  

%Note that, unlike \cite{KK}, we allow more than 
%There is at most one  morphism $f=f_{c_2,c_1}\in \mor(\CC)$ so that $f(c_1)\subset c_2$. 
%We say that a {\em hyperbolic polyhedral complex}  is {\em regular} if it satisfies: 

Axiom 2. For every $c_1, c_2\in \Ob(\CC)$ there exists at most one morphism $f=f_{c_2,c_1}\in \mor(\CC)$ so that $f(c_1)\subset c_2$. 
\end{definition}

Objects of a polyhedral complex $\CC$ are called {\em faces} of $\CC$ and the morphisms of $\CC$ are called {\em incidence maps} of $\CC$. 
A {\em facet} of $\CC$  is a face $P$ of $\CC$ so that 
for every morphism $f: P\to Q$ in $\CC$, $f(P)=Q$.  A {\em vertex} of $\CC$ is a zero-dimensional face. The {\em dimension} 
$\dim(\CC)$ of $\CC$ is the supremum of dimensions of faces of $\CC$. 
A polyhedral complex $\CC$ is called {\em pure} if the dimension function is constant on the set of facets of $\CC$; 
the constant value in this case is the dimension of $\CC$. A {\em subcomplex} of  $\CC$  is a full subcategory of $\CC$.   
If $c$ is a face of a complex $\CC$ then $\Resid_{\CC}(c)$, the {\em residue} of $c$ in $\CC$, is the minimal 
subcomplex of $\CC$ containing all faces  $c'$ such that there exists an incidence map $c\to c'$. For instance, if $c$ is a vertex of 
$\CC$ then its residue is the same as the {\em star} of $c$ in $\CC$; however, in general these are different concepts. 

\medskip 
We generate the equivalence relation $\sim$ on a polyhedral complex $\CC$ by declaring that $c\sim f(c)$, where $c\in \Ob(\CC)$ and 
$f\in \mor(\CC)$. This equivalence relation also induces the equivalence relation $\sim$ on points of faces of $\CC$.

If $\CC$ is a polyhedral complex, its {\em poset} $\operatorname{Pos}(\CC)$ is the partially ordered set $\Ob(\CC)$ with the relation $c_1\le c_2$ iff $c_1\sim c_0$ so that $\exists f\in \mor(\CC), f: c_0\to c_2$. 

We define the {\em underlying space} or {\em amalgamation} $|\CC|$ of a polyhedral complex $\CC$ as the topological space which is obtained from the disjoint union 
$$
\amalg_{c\in \Ob(\CC)}\ c
$$
by identifying points using the equivalence relation: $\sim$. We equip  $|\CC|$ with the quotient topology.

\begin{definition}\label{difference}
If $\CC$ is a polyhedral complex and $\CB$ is its subcomplex. For $c\in \Ob(\CC)$ define the polyhedron 
$$
 c^{\CB}:=c\setminus \bigcup_{b\le c, b\in \CB}  f(b), \hbox{~~where~~}  f: b\to c, f\in \mor(\CC). 
$$
For a morphism $f\in \mor(\CC)$, $f: c_1\to c_2$, we set $f^{\CB}: c_1^{\CB}\to c_2^{\CB}$ be the restriction of $f$. 
We define the {\em difference complex} $\CC - \CB$ as the following polyhedral complex:
$$
\Ob(\CC - \CB)=\{ c^{\CB}: c \in \Ob(\CC)\},   
$$
$$
\mor(\CC - \CB)= \{f^{\CB}: c_1^{\CB}\to c_2^{\CB},  \hbox{where~} f\in \mor(\CC), f: c_1\to c_2\}. 
$$
\end{definition}

A complex $\CC$ is said to be {\em finite} it has only finitely many objects and morphisms. A complex $\CC$ is {\em locally finite} 
if for every face $a\in Ob(\CC)$ the sets of morphisms
$$
\{ f: a\to b, b\in Ob(\CC)\}%, \quad \{ f: b\to a\}, %\quad b\in Ob(\CC)
$$ 
is finite. The key example of a hyperbolic polyhedral complex used in this paper (the Dirichlet tiling of $\H^n$) 
will be infinite but locally finite. In this paper we will be exclusively interested in locally finite complexes.

\begin{definition}
Let $\CC$ be a pure $n$-dimensional polyhedral complex. 
The {\em nerve} $\Nerve(\CC)$ of $\CC$ is the simplicial complex whose vertices are facets of $\CC$ 
(the notation is $v=c^*$, where $c$ is a facet of $\CC$); 
distinct vertices $v_0=c_0^*,...,v_k=c_k^*$ or $\Nerve(\CC)$ span a $k$-simplex if there exists an $n-k$-face $c$ of $\CC$ and 
incidence maps $c\to c_i, i=0,...,k$. The simplex $\sigma=[v_0,...,v_k]$ then is said to be {\em dual} to the face $c$.    
\end{definition}

 Similarly to \cite{KK} we have: 

\begin{lemma}
If $\CC$ is locally finite then $|\CC|$ is homotopy-equivalent to $|\Nerve(\CC)|$. 
\end{lemma}

\begin{definition}
A hyperbolic polyhedral complex $\CC$ is {\em simple} if: 

(1) $\CC$ is pure and $\dim(\CC)=n$,  

(2) For $k=0,...,n$ and every $k$-face $c$ of $\CC$, 
$\Nerve(\Resid_{\CC}(c))$ is isomorphic to the complex $\CC(\Delta^{n-k})$. 
\end{definition}

For a polyhedral complex $\CC$ we define its {\em $k$-skeleton}, to be the subcomplex $\CC^{(k)}$ consisting of faces of dimension $\le k$. 
For a pure $n$-dimensional complex $\CC$ we define its {\em derived complex} complex $\CC'$ by: $\CC':= \CC - \CC^{(n-3)}$.  
We say that $\CC$ is {\em weakly simple} if the derived complex $\CC'$ is simple.  
In other words, every $n-2$-dimensional face is incident to exactly 3 facets and every $n-1$-dimensional face 
 is incident to exactly 2 facets.

The point of considering derived complexes is that 
if $|\CC|$ is a manifold at every point of $\CC^{(n-3)}$, then 
$$
\pi_1(|\CC|)\cong \pi_1(|\CC'|). 
$$
Thus, in this situation, passing to the derived complex does not change the fundamental group, while proving simplicity for the derived complex is much easier than for the original one. 

\medskip
{\bf Voronoi tiling of $\H^n$.} 

\begin{defn}
Let $Y\subset \H^n$ be a locally finite subset (i.e., every compact in $\H^n$ contains only finitely many points of $Y$). 
The {\em Voronoi tiling} $\CV(Y)$ of $\H^n$ associated with $Y$ is defined by:  For each $y\in Y$ take the {\em Voronoi cell} 
$$V(y):=\{x\in \H^n: d(x,y)\le d(x,y'), \forall y'\in Y\}.$$ 
Thus, each cell $V(y)$ is given by the collection of non-strict linear inequalities $d(x,y)\le d(x,y')$. 
Then each cell $V(y)$ is a closed (possibly unbounded) polyhedron in $\H^n$. 
The union of Voronoi cells is the entire $\H^n$. Assuming that each $V(y)$ has only finitely many faces, we thus obtain the polyhedral complex, called the {\em Voronoi complex}, $\CV(Y)$ using the polyhedra $V(y)$ as facets and faces of facets as faces of $\CV(Y)$. 
\end{defn}

A special case of this construction is given by orbits of a discrete convex-cocompact subgroup $\Ga< PO(n,1)$: If $x\in \H^n$ is a point 
not fixed by any $\ga\in \Ga\setminus \{1\}$, then the Dirichlet tiling $\CD_x(\Ga)$ is the same as 
Voronoi tiling with respect to the set $Y:=\Ga\cdot x$. We will use the same notation $\CD_x(\Ga)$ for the associated hyperbolic 
polyhedral complex, the {\em Dirichlet complex}.  (Recall that, since $\Ga$ is convex-cocompact, every  
$D_{\ga x}$ is a convex hyperbolic polyhedron in our sense since it has only finitely many faces.)

\medskip
We note that if $c$ is a face of $D_x=D_x(\Ga)$, then its stabilizer in $\Ga$ has to be finite, otherwise $D_x$ would fail to be a fundamental domain. 
In particular, the stabilizer  of $c$ consists entirely of elliptic elements.

\begin{lemma}\label{lem:regular}
Let $\CD_x(\Ga)$ be the Dirichlet complex of a convex-cocompact group $\Ga< PO(n,1)$. Then $\Ga$ contains a finite-index 
torsion-free normal subgroup $\Ga'$  so that $\left(\CD_x(\Ga)\right)/\Ga'$ is a   
hyperbolic polyhedral complex. 
\end{lemma}
\proof First, since $\Ga$ is convex-cocompact, it is also finitely-generated, see e.g. \cite{Bowditch}. Hence, by Selberg's Lemma \cite{Selberg}, 
$\Ga$ contains a torsion-free subgroup $\Ga_1$ of finite index. One could  now take the quotient 
$\left(\CD_x(\Ga)\right)/\Ga_1$: It satisfies all properties of a polyhedral complex, except Axiom 2 could fail: If $\tilde{c}_1\le \tilde{c}_2$ are incident 
faces of $\CD_x(\Ga)$, we could have some $\ga \in \Ga_1\setminus \{1\}$ such that $\ga(\tilde{c}_1)\le \tilde{c}_2$. Dividing by $\Ga_1$ we then 
would have more than one morphism $c_1\to c_2$, where $c_i$ is the projection of $\tilde{c}_i, i=1,2$. We will see below how to eliminate  
such elements $\ga$ by passing to a further finite index subgroup in $\Ga_1$. 

Since $D_x$ has only finitely many faces, there are only finitely many nontrivial 
elements $\ga_i\in \Ga_1, i=1,...,m,$ so that  $\ga_i D_x\cap D_x\ne \emptyset$. Since $\Ga_1$ is residually finite, it contains a finite-index 
subgroup $\Ga_2$ so that $\ga_i\notin \Ga_2, i=1,...,m$. Lastly, we take $\Ga'<\Ga_2$ a finite index subgroup which is normal in $\Ga$. 
Then for each $\ga\in \Ga'\setminus \{1\}$, $\ga D_x\cap D_x=\emptyset$. By normality of $\Ga'$ in $\Ga$, we also have 
$$
\ga D_{\al x}\cap D_{\al x}=\emptyset
$$
for all $\al\in \Ga$. This implies that Axiom 2 holds for the quotient complex $\left(\CD_x(\Ga)\right)/\Ga'$. \qed

\medskip
{\bf Weak simplicity criterion for Dirichlet complexes.} 

\begin{lemma}\label{L2.1}
The Dirichlet tiling $\CD_x(\Ga)$ is weakly simple provided that for every $y\in \partial D_x\subset \H^n$ and every collection of elements 
$\ga_1,...,\ga_k\in \Ga$ so that 
$$
\dim\left( \bigcap_{i=1}^k Bis(x, \ga_i x) \right) =n-2, 
$$
and 
$$
y\in \bigcap_{i=1}^k Bis(x, \ga_i x),  
$$
we have $k=2$. 
\end{lemma}
\proof Since $\Ga$ acts transitively on the facets of the tiling $\CD_x(\Ga)$, it suffices to prove weak simplicity of $\CD_x(\Ga)$ 
along codimension 2 cells $E$ contained in the boundary of $D_x$. Let $\ga_0:=1, \ga_1,...,\ga_k\in \Ga$ be the elements of 
$\Ga$ such that $E$ is contained in $\ga_i(D_x)$, $i=0,...,k$. Weak simplicity of $\CD_x$ then means that $k=2$. 
We relabel the elements $\ga_i$ above so that 
$$
\ga_{i+1} D_x\cap \ga_i D_x
$$
is a codimension 1 face $F_i$ for $i=0,...,k$, where $i$ is taken modulo $k$. Then $F_i$ is contained in the bisector
$Bis(\ga_i x, \ga_{i+1}x)$, $i=0,...,k$. Therefore, for every $y\in E$, 
$$
d(y, \ga_i x)= d(y, x),
$$
that is,
$$
y\in \bigcap_{i=1}^k Bis(x, \ga_i x)
$$
and, hence,
$$
E\subset \bigcap_{i=1}^k Bis(x, \ga_i x). \qed 
$$

\begin{remark}
The same proof, of course, yields the simplicity criterion for $\CD_x(\Ga)$: It is simple if and only if 
for every $y\in \partial D_x$, the bisectors $Bis(x, \ga_i(x))$ passing through $y$ are transversal. 
% and every collection of bisectors $Bis(x, \ga_i(x))$ passing through $y$, these bisectors are transversal. 
\end{remark}

{\bf Linear algebra problems.} Let $F$ be a field.  
%and $G< GL(n, F)$ be a subgroup (in the most interesting cases, $G$ is an algebraic 
%subgroup of $GL(n,F)$). 
For a subset $\ul{A}=\{A_1,...,A_k\}\subset Mat_{n,n}(F)$, $k\le n$, we define the  map
\begin{equation}\label{eq:B}
B=B_{\ul{A}}: F^n \to Mat_{n,k}(F)
\end{equation}
by
$$
x\mapsto (B_1 x,...,B_k x)
$$
where $B_i=A_i-I, i=1,...,k$ and we regard vectors $B_i x$ as columns of the matrix $B_{\ul{A}}$. We say that the map $B$ 
and the set $\ul{A}$ are {\em singular} if for every $x\in F^n$, $\rank(B(x))<k$. We note that the image of the map $B$ is a 
linear subspace of $Mat_{n,k}(F)$. The problem of describing linear subspaces of $Mat_{n,k}(F)$ consisting of matrices 
of rank $<k\le n$ has a long history, see \cite{Lovasz} for a survey.  

If we do not make any restrictions on the matrices $A_i$, then the problem of describing singular 
$k$-tuples $\ul{A}$  is essentially equivalent to the problem of describing linear subspaces of $Mat_{n,k}(F)$ 
and is hopelessly complicated. Suppose, however, one takes $A_i$ from an algebraic subgroup $G< GL(n,F)$, e.g., $G=O(n,F)$. 

\begin{problem}\label{prob:LA}
Let $F=\R$. Describe singular $k$-tuples $\ul{A}$ of matrices $A_i\in O(n,1)^\uparrow$. In particular, suppose that no matrices in $A_i$ 
share a common eigenvector. Is it true that in this case  $\ul{A}$ is nonsingular? 
\end{problem}

A positive answer would be a key step towards proving 

\begin{conjecture}
Suppose that $\Ga< PO(n,1)$ is a discrete subgroup of the class $\CK$. Then for generic $x\in \H^n$ the Dirichlet complex $\CD_x(\Ga)$  is simple. 
\end{conjecture}

Note that the linear maps $B=B_{\ul{A}}: F^n\to Mat_{n,k}(F)$  are injective provided that the linear 
transformations $A_i$ do not have a common fixed vector. In this case, Problem \ref{prob:LA} becomes a special case of the problem of describing 
$k$-dimensional linear subspaces of $Mat_{n,k}(F)$ (with $k\le n$) consisting of matrices of rank $<k$. {\em Rank} of such a subspace 
is the maximal rank of a matrix which belongs to this subspace.  

Linear subspaces of $Mat_{n,k}(F)$ of rank $1$ are easy to describe. Classification of subspaces of ranks $2$ and $3$ was given in 
\cite{Atkinson, EH}. It is easy to see that the classes of 
{\em primitive} subspaces of rank $\le 3$ described in \cite{Atkinson, EH} (with $F=\C$) 
do not appear as images of maps $B_{\ul{A}}$, where $\ul{A}=\{A_1,A_2,A_3,A_4\}$ are in $O(4,\C)$. However, 
it is unclear how to deal with the non-primitive subspaces. For instance, it is unclear if there are (pairwise noncommuting) 
elements $A_i$ of $O(4,\C)$ so that the matrices $B_i$ are linearly dependent as elements of $Mat_{4,4}(\C)$. Such quadruples 
would correspond to the case when
$$
\bigcap_{x\in F^4} Ker( B(x))\ne 0. 
$$

\section{Complexes of varieties}\label{sec:varieties}

Our discussion here closely follows \cite{KK}. 

\begin{defn} Let $\VV$ denote either the category of varieties
(over a fixed field $k$) or the category of topological spaces.

Let $\CC$ be a finite hyperbolic polyhedral complex. 
A {\em $\VV$-complex based on $\CC$} is a functor $\Phi$ from $\CC$ to $\VV$ so that 
 morphisms $c_i\to c_j$ go to closed embeddings 
$\phi_{ij}:\Phi(c_i)\to \Phi(c_j)$. By abuse of terminology, we will sometimes refer to the image category $im(\Phi)$ as a $\VV$-complex based on $\CC$. 
We will use the notation $X_i$ for $\Phi(c_i)$. The varieties $X_i$ will be called {\em strata} of the complex of varieties $im(\Phi)$. 

We call the functor $\Phi$ {\it strictly faithful} if the following holds:

 If $x_i\in \Phi(c_i)$, $x_j\in \Phi(c_j)$ and $\phi_{ik}(x_i)=\phi_{jk}(x_j)$ 
for some $k$ then there is an $\ell$
and $x_{\ell}\in \Phi(c_{\ell})$ such that
$\phi_{\ell i}(x_{\ell})=x_i$ and $\phi_{\ell j}(x_{\ell})=x_j$.

The relation $x_i\sim \phi_{ij}(x_i)$ for every $i, j$ and $ x_i\in X_i$
generates an equivalence relation
on the points of $\amalg_{i\in I} \Phi(c_i)$, also denoted by $\sim$. 
\end{defn}

In the category of topological spaces, the direct limit (or push-out) $\lim\Phi(\CC)$
of the diagram $\Phi(\CC)$  exists and its points are identified with
 $\bigl(\amalg_{i\in I} \Phi(c_i)\bigr)/\sim $.

For example, suppose that $\Phi_{taut}$ is the tautological functor which identifies each face of $\CC$ 
with the corresponding underlying topological space. Then $\lim \Phi_{taut}(\CC)$ is nothing but $|\CC|$. 

As in \cite{KK} we have:

\begin{lemma}\label{pi_1-cong}
Suppose that $\Phi$ is strictly faithful and $\Phi(\CC)$ consists of cell complexes and cellular maps of such complexes. Then 
 $\pi_1\bigl(\lim\Phi(\CC)\bigr)\cong \pi_1(|\CC|)$ provided that each $\Phi(c), c\in \Ob(\CC)$ is $1$-connected. \qed
\end{lemma}

The following result was proven in \cite{KK}, Proposition 31 for Euclidean polyhedral complexes, but the same proof applies to hyperbolic 
complexes: 

\begin{proposition}\label{prop:variety}
Let $\Phi: \CC \to \VV$ be a complex of varieties based on a finite hyperbolic polyhedral complex. 
Assume that for each $k$ and each $J\subset I$ the subvariety  $\cup_{j\in J}\im(\phi_{jk})\subset X_k$ is seminormal. For instance, this assumption holds if this union is a divisor with normal crossings. Then the direct limit $\lim\Phi(\CC)$ exists in the category of varieties. Furthermore, as a topological space, 
$\lim\Phi(\CC)$ is homeomorphic to the topological push-out $\lim\Phi_{top}(\CC)$. 
\end{proposition}

{\bf Complex of projective spaces.} Let $\CC$ be a finite simple hyperbolic polyhedral complex; we set $C:=|\CC|$. As in \cite{KK}, we define the functor 
$\Phi: \CC\to \VV$ which sends each $c$ to the projective space $\CP_{c}:=\Span(c)\times \{c\}$ and each morphism $c_i\to c_j$ 
to the linear map of the projective spaces $F_{c_j,c_i}: \CP_{c_i}\to \CP_{c_j}$ which restricts to the original morphism   $c_i\to c_j$. 
We let $\CP=\CP_{\CC}$ denote $\im(\Phi)$. 

The following definition is taken from \cite{KK}: 

\begin{defn}[Parasitic subspaces]\label{paras.int.defn}
Let  $\sigma:=(c_1, c_2,...,c_k)$ be a tuple of faces
incident to a face $c$. Consider the intersections 
$$
I_{c, \sigma}:= \cap_{i=1}^k F_{c,c_i}(\CP_{c_i})\subset \CP_{c}
$$
such that there is no face $c_0$ such that
$I_{c, \sigma}=F_{c,c_0}(\CP_{c_0})$ and $c_0$ is incident to all the 
$c_1, c_2,...,c_k$. Then the subspace $I_{c, \sigma}\subset \CP_{c}$ 
is called a {\em parasitic intersection} in $\CP_{c}$.

Note, however, that this collection of parasitic intersections in spaces $\CP_c, c\in Ob(\CC)$, 
is not stable under applying morphisms $F_{c',c}$ and taking preimages under these morphisms. 
We thus have to {\em saturate} the collection of parasitic intersections using the morphisms $F_{c,c'}$. 
This is done as follows. Let $T$ denote the push-out of the category $im(\CP_{top})$, where we regard each $\CP_b, b\in Ob(\CC)$,   
as a topological space, so the push-out exists. Then for each $a\in Ob(\CC)$ we have the (injective) projection map 
$\rho_a: \CP_a\to T$.  For each parasitic intersection $I_{c,\si}\subset \CP_c$, we define
$$
I_{c,\si,a}:=\rho_a^{-1} \rho_c(I_{c,\si}). 
$$
We call such $I_{c,\si,a}$ a {\em primary parasitic subspace} in $\CP_a$. It is immediate that each primary parasitic subspace in $\CP_a$ 
is a projective space linearly embedded in $\CP_a$. With this definition, the collection of parasitic subspaces $I_{c,\si,a}$ 
is stable under taking images and preimages of the morphisms $F_{c,c'}$. 
\end{defn}

For the purposes of this paper, we will need an equivariant version of the above definition. %Suppose that $\CC$ is finite. 
Each $\CP_c$ embeds in $|\CP|$, the push-out of $\Phi_{top}(\CP)$.

by abusing the notation we retain the notation $\CP_c$ for the image in $|\CP|$. 
Let $\Theta$ be a group acting faithfully and isometrically on $\CC$. This action extends to a faithful (linear) action $\Theta \acts \CP$. For 
$\theta\in  \Theta\setminus \{1\}$, consider the fixed-point set $Fix(\theta)$ of $\theta$ in $|\CP|$. For each point $p\in Fix(\theta)$ 
we take the smallest face $c=c(p)$ such that $p\in \CP_{c}$. By minimality, $\theta(c)=c$. We then let 
$Fix_{c}(\theta):=Fix(\theta)\cap \CP_{c}$, it is a finite union of disjoint projective subspaces in $\CP_{c}$. 
We obtain a  set of projective spaces $\CP_{c_i}$ (in the example we are mostly interested in this is a single projective space) 
so that
$$
\forall i, \quad \theta(c_i)=c_i, \quad \hbox{and}\quad Fix(\theta)\subset \bigcup_{i} Fix_{c_i}(\theta). 
$$

\begin{assumption}\label{asm:cartan}
For every $\theta\in \Theta\setminus \{1\}$ which does not act freely on $C$, $\theta^2=1$ and, moreover, for every $c_i$ as above, 
$Fix_{c_i}(\theta)= \{p_i\} \sqcup p_i^\perp$, where $p_i^\perp$ 
is a codimension 1 projective subspace in $\CP_{c_i}$ so that 
$$
c_i\cap p_i^\perp=\emptyset. 
$$ 
Furthermore, each $p_i$ belongs to exactly three faces of $\CC$: Two facets $a_i, b_i$ and one $\theta$-invariant 
codimension 1 face $c_i$ incident to $a_i, b_i$. Note that we do not assume that $\CC$ is finite here. 
\end{assumption}

%Even though it is not needed for our paper, we generalize this setting as follows. 
%\begin{definition}
%We say that the action $\Theta \acts \CC$ is {\em fine} if for every $\theta\in \Theta$ and every $c\in Ob(\CC)$, $Fix_c(\theta)$ is a finite 
%union of pairwise disjoint projective subspaces $P_i$  so that each $P_i$ is either disjoint from $c$ 
%\end{definition}

 Thus, $p_i^\perp\subset \CP_{c_i}=\CP_{a_i}\cap \CP_{b_i}$. However, a priori, 
$p_i^\perp\subset |\CC|$ could intersect other strata as well. We would like to eliminate these intersections. 
 (Notice that we will be ignoring parasitic subspaces of $\CP_{c_i}$ which could cross $p_i^\perp$.) 
Namely, for each  face $c_i$ we consider strata $\CP_{e_i}\subset  \CP_{c_i}$, %of dimension $<\dim(c_i)$, 
where $e_i$'s are  
proper faces of $c_i$. For each $e_i$ define $Q_{e_i,p_i}:=\CP_{e_i}\cap p_i^\perp$ provided that this intersection is nonempty. 
We would like to get rid of the subspaces $Q_{e_i,p_i}$, so we declare the intersections  
$Q_{e_i,p_i}\in \CP_{c_i}$ to be {\em secondary parasitic subspaces}. We saturate the collection of secondary parasitic subspaces as we 
did in the primary case.   
 
\begin{remark}
The case we are mostly interested in is when $\dim(\CC)=3$, so each $e_i$ is an edge and each secondary parasitic subspace 
$Q_{e_i,p_i}$ is just a point. 
\end{remark}

We now assume that $\CC$ and $\Theta$ are both finite. We proceed as in \cite{KK} and blow up all the parasitic subspaces: 
We first blow up all primary parasitic subspaces (by induction on dimension) and then blow up all secondary parasitic 
subspaces (again, by induction on dimension). 
The construction is canonical, so the group $\Theta$ 
continues to act on the blow-up $b\CP$. By applying Proposition \ref{prop:variety}, we conclude that the 
$\Theta$-equivariant push-out $X:=|b\CP|$ exists in the category of projective varieties and is equivariantly 
homeomorphic to the topological push-out. The same arguments as in \cite[\S 5]{KK} show that the variety $X$ is projective. 
As in \cite{KK}, the variety $X$ has only normal crossing singularities. 

Furthermore, by the construction, in view of the Assumption \ref{asm:cartan}:

\begin{lemma}\label{lem:nonfree}
1. $\theta\in \Theta$ has a fixed point in $C$ if and only if $\theta$ has a fixed point in $X$. Such $\theta$ has order 2. 

2. For every $\theta\in \Theta\setminus \{1\}$ which does not act freely on $X$, every 
component of $Fix(\theta)\subset X$ is contained in the intersection of exactly two top-dimensional strata (of dimension $n$) 
intersecting normally. The local models for the action of $\theta$ are described below.
\end{lemma}

Let $y_1,...,y_{n+1}$ be coordinates on $\C^{n+1}$. Then:  

1. Near an isolated fixed point $p_i$:
$$
y_1  y_2=0, \quad \theta(y_1,y_2,...,y_{n+1})=  (y_2,y_1, -y_3, -y_4...,-y_{n+1}). 
$$

2. Along the $n-2$-dimensional component $bFix_{c}(\theta)$, (the blow-up of $Fix_{c}(\theta)$), where $\dim(c)=n-1$ :
$$
y_1  y_2=0, \quad \theta(y_1,y_2, y_3,...,y_{n+1})=  (y_2,y_1, -y_3, y_4...,y_n , y_{n+1}). 
$$  

The case we are mostly interested in is when $n=3$, so the latter action becomes:
$$
y_1  y_2=0, \quad \theta(y_1,y_2,y_3,y_{4})=  (y_2,y_1, -y_3, y_4). 
$$

We will refer to these singularities together with the $\Z_2$-actions as $(Y_1,0)$ and $(Y_2, 0)$ respectively. Notice that 
if we blow up the origin in $(Y_1,0)$, then we obtain singularity of the 2nd type. 

Notice that $Y_2$ splits equivariantly as the product $Y\times \C$, where 
$$
Y=\{ (y_1,y_2,y_3)\in \C^3: y_1 y_2=0\},  \theta(y_1,y_2,y_3)=  (y_2,y_1, -y_3)
$$
and the action of $\theta$ on the remaining factor $\C$ is by the identity. Hence, $Y_2/\Z_2\cong Y/\Z_2 \times \C$. The variety $Y/\Z_2$ is 
a normal crossing along the line $y_1=y_2=0$ away from the origin. I am grateful to J\'anos Koll\'ar for providing the proof of the following: 

\begin{lemma}
The germ of $Y/\Z_2$ at the origin is isomorphic to the Whitney umbrella $u^2=w v^2$.  
\end{lemma}
\proof The ring of invariants $\C[y_1,y_2,y_3]^{\Theta}$ is generated by the polynomials 
$y_1y_2, y_1+y_2, y_3(y_1-y_2), y_3^2$ subject to the equation
$$
\left(y_3(y_1-y_2)\right)^2=y_3^2\left((y_1+y_2)^2-4y_1y_2\right). 
$$
Dividing this ring by the ideal generated by $y_1y_2$ we obtain the ring $Q$ with the generators 
$u:=y_1+y_2, v:=y_3(y_1-y_2), w:=y_3^2$ subject to 
the equation
$$
\left(y_3(y_1-y_2)\right)^2=y_3^2 (y_1+y_2)^2.   
$$
Equivalently, $Q$ is generated by $u,v,w$ subject to the equation
$$
v^2=w u^2.   
$$
However, this is the quotient ring of the Whitney umbrella.  \qed

\medskip
As in \cite{KK}, the blow-up $b\CP$ is strictly faithful and, hence, $\pi_1(X)\cong \pi_1(C)$, where $C=|\CC|$. We let $N$ denote the group 
$\pi_1(X)= \pi_1(C)$. 

\begin{proposition}\label{prop:centerless}
Suppose that $N$ has trivial center. Then $\pi_1(X/\Theta)\cong \pi_1(C/\Theta)$. 
\end{proposition}
\proof We have a $\Theta$-equivariant isomorphism of fundamental groups $\pi_1(C)\to \pi_1(X)$.  Thus, considering the quotient--orbihedra 
$\CO_C:=C/\Theta$ and $\CO_X:=X/\Theta$, we obtain group extensions
$$
1\to N\to \pi_1(\CO_C)\to \Theta\to 1, \quad 1\to N\to \pi_1(\CO_X)\to \Theta\to 1 
$$
where the homomorphisms $\psi_i: \Theta\to Out(N)$ associated with the actions of $\Theta$ on $\pi_1(C)$ and $\pi_1(X)$ are the same. 
Since $N$ has trivial center, by Corollary 6.8 in \cite[Chapter IV]{Brown}, the group extensions above are naturally isomorphic. 
Define normal subgroups $F_C, F_X$ of the groups $\pi_1(\CO_C), \pi_1(\CO_X)$ to be the normal closures of 
the elements of the respective groups, which do not act freely on the universal covers of the orbihedra $\CO_C$ and $\CO_X$. 
By Armstrong's theorem \cite{Armstrong}, the fundamental groups of $C$ and $X$ are obtained from the orbihedral fundamental groups 
$\pi_1(\CO_C), \pi_1(\CO_X)$ by dividing by the subgroups $F_C, F_X$. We claim that 
the isomorphism $\pi_1(\CO_C)\to \pi_1(\CO_X)$ carries $F_C$ to $F_X$ isomorphically. 

Indeed, let $\tilde{C}\to C$ and $\tilde{X}\to X$ denote the universal covers of $C$ and $X$ respectively. The space $\tilde{C}$ 
has a natural structure push-out of a hyperbolic polyhedral complex $\tilde{\CC}$, while  
$\tilde{X}$ has a natural structure of push-out of a complex of varieties $\tilde\CP$ based on $\tilde{\CC}$. The strata of $\tilde{X}$ project 
isomorphically to the strata of $X$ since the latter are simply-connected. 

Suppose now that, say, $\tilde\theta\in \pi_1(\CO_C)\setminus \{1\}$ is a lift of $\theta\in \Theta$ 
has a fixed point $\tilde{p}$ in the universal cover of $C$. The isomorphism $\pi_1(C)\to \pi_1(X)$ is induced by the natural embedding of the universal covers $\iota: \tilde{C}\to \tilde{X}$. Therefore, such $\tilde\theta$ also fixes the point $\iota(p)\in X$. Conversely, if  
$\tilde\theta\in \pi_1(\CO_X)\setminus \{1\}$  fixes a point $\tilde{q}$ in the universal cover of $X$, then 
$\tilde{q}$ belongs to a minimal stratum $\tilde{X}_i$ of $\tilde{X}$, which corresponds to a face $\tilde{c}_i$ of $\tilde{\CC}$. Then $\tilde\theta$ has to 
preserve $\tilde{X}_i$ and, hence, $\tilde{c}_i$. The projection $\tilde{X}_i\to X_i\subset X$ is an isomorphism conjugating the action of 
$\tilde\theta$ to the action of $\theta$. Since $\theta$ was fixing a point in $c_i\in Ob(\CC)$ (where $c_i$ is the image of $\tilde{c}_i$ under 
the projection $\tilde{C}\to C$), we conclude that $\theta$ also fixes a point in $\tilde{c}_i$. Proposition follows. \qed 

\medskip 
{\bf Dimension reduction.} Let $V$ be the variety obtained from $X/\Theta$ as follows: We first equivariantly blow up in $X$ all isolated fixed 
points of involutions $\theta\in \Theta$ and then divide the resulting variety by $\Theta$. The quotient has only normal crossing singularities and singularities of the 2nd type, more precisely, of the type $Y_2/\Z_2$. These singularities split as the product $Y/\Z_2 \times \C$, where $Y/\Z_2$ is a Whitney umbrella. We now embed $V$ in the projective space and intersect it with a generic hyperplane. The result is a projective surface 
$V$ whose singularities 
are only normal crossings and Whitney umbrellas. Furthermore, by Lefschetz Hyperplane section theorem $\pi_1(V)\cong \pi_1(W)$, 
see \cite[p. 27]{gm-book}. Since $V$ was irreducible, so is $W$.

\section{Generic transversality of triples of bisectors in $\H^n$ }%Proof of Theorem \ref{thm:nonsingular}} 

The main result of this section is

\begin{theorem}\label{thm:simple}
Let $A_1, A_2, A_3\in O(n,1)^\uparrow$ be distinct nontrivial elements of a group $\Ga< O(n,1)$ of the class $\CK$, $n\ge 2$. 
Assume also that $A_1, A_2, A_3$ do not generate a cyclic group. Then for generic $x\in H$, the vectors 
$$
B_i(x)=A_i(x)-x, i=1,2,3
$$
are linearly independent. 
\end{theorem}
\proof Linear dependence of the vectors $B_i(x)$ is equivalent to the condition that $\rank(B_{\ul{A}}(x))\le 2$, which, in turn, 
is expressed in terms of vanishing of determinants of $3\times 3$ minors of the 
$(n+1)\times 3$ matrix $B_{\ul{A}}(x)$. 
Therefore, the set of $x\in H$ such that $\rank(B_{\ul{A}}(x))\le 2$ is an algebraic subset. Hence, it is either the entire 
$H$ or it is a closed set with empty interior. 

We suppose therefore, that for every $x\in H$ the vectors $B_i(x), i=1,2,3$ 
are linearly dependent. Then, by linearity, the same is true for all $x\in \CL^\uparrow$.   
Since $\CL^\uparrow$ is Zariski dense in $\R^{n+1}$, the same conclusion holds for all $x\in \R^{n+1}$.

We let $\Om\subset C^+\subset \R^{n+1}$ 
denote the set of $x\in C^+$ such that $\rank(B_{\ul{A}}(x))= 2$. 
Our first goal is to understand the complement of $\Om$, i.e., the set of $x\in C^+$ such that all the three vectors $B_i(x)$ are multiples of each other. 
We will consider a (seemingly) larger set 
$$
\Si=\Si_{12}\cup \Si_{23}\cup \Si_{31}\subset C^+$$
where $\Si_{ij}= \{x\in C^+: \dim(Span(B_i(x), B_j(x)))\le 1\}$.

\begin{lem}\label{L3.1}
Let $x\in C^+$ be a nonzero vector. Then $x\in \Si_{12}$ iff one of the following holds:

1. $A_1 x=x$ or $A_2 x=x$ or %(i.e., $B_i(x)=0$ for $i=1$ or $i=2$). 
$A_1 x=A_2 x$, i.e., $x$ is fixed by $A_2^{-1} A_1$.  This can happen only if $x\notin C$. 

2. $x$ is a common eigenvector of $A_1, A_2$. This can happen only if $x\in C$.%, provided $x\cdot x=0$ then 
\end{lem}
\proof If $B_i(x)=0$ then $A_i(x)=x$. We, thus, will assume that $B_i(x)\ne 0, i=1,2$. 

The condition $x\in \Si_{12}$ then is equivalent to 
$$
B_1(x)=\mu B_2(x), \mu\ne 0. 
$$
In other words, 
$$
%A_1 x - x= \mu A_2 x -\mu x, \quad 
A_1 x - \mu A_2 x= (1-\mu) x. 
$$

1. If $\mu=1$ then $A_1x=A_2 x$, $A_2^{-1} A_1 x=x$. Furthermore, every $x$ satisfying these properties belongs to $\Si_{12}$ (by taking $\mu=1$). 

2. Suppose now that $\mu\ne 1$. 

a. If  $x\cdot x\ne 0$, Lemma \ref{L2.0} then implies that $A_ix= x, i=1,2$, contradicting our assumption $B_i(x)\ne 0, i=1,2$. 

b. Suppose that $x\cdot x=0$. Then linear dependence of the vectors $x, A_1 x, A_2 x$ (which belong to the conic $C$) 
implies that they belong to a common line in $\R^{n+1}$. In particular, $x$ is a common eigenvector of $A_1, A_2$. \qed

\begin{cor}\label{C1}
The set $\Si_{ij}$ is a finite union of lines. %linear subspaces of dimension $\le 2$. 
In particular, $\Si$ does not separate $C^+$ and, thus, the open set  $\Om$ is connected. 
\end{cor}
\proof We need to observe two things: First, $\Ga$ does not contain elliptic elements besides Cartan involutions. Hence, 
fixed-point sets and eigenspaces of $A_1, A_2, A_2^{-1}A_1$ in $C^+$ are at most lines. 
Secondly, since $n\ge 2$, no line can separate $C$. Now, the statement follows from Lemma \ref{L3.1}. \qed 

\medskip
Recall that we are assuming that for all $x\in C^+$, the vectors $B_i(x), i=1,2,3$ are linearly dependent. Thus, there exist (possibly multivalued and discontinuous) functions $\al_i(x), i=1,2,3$ so that for all $x\in C^+$
\begin{equation}\label{eq:ld}
\sum_{i=1}^3 \al_i(x) B_i(x)=0 . 
\end{equation}
If for $x\in C^+$ one can take $\al_k(x)=0$ then $x\in \Si_{ij}$, $\{i,j,k\}=\{1,2,3\}$. In particular, for each $x\in \Om$ all the quantities $\al_i(x)$ are nonzero. 
Hence, we can select (say, by setting $\al_1(x)\equiv 1$) 
nonvanishing continuous functions $\al_i(x), i=1,2,3$, $x\in \Om$, so that \eqref{eq:ld} holds. We fix these 
functions from now on. 

%\begin{cor}\label{C2}
%Away from an (at most) 2-dimensional subset $\Si\subset \R^4$, the vectors $B_i(x), i=1,2,3$ are all nonzero 
%and span a 2-dimensional subspace. Thus, there exists a well-defined (up to scaling) triple of continuous functions $\al_i(x), i=1,2,3$ on $\Om$,  so that 
%$$\sum_{i=1}^3 \al_i(x) B_i(x)=0, \quad \forall x\in \Om. $$
%Furthermore, the zero level sets of the functions $\al_i(x)$ do not separate $\R^4$ and $C$. \end{cor}

%Consider the map $\Phi_{A_1,A_2}: x\mapsto \star(B_1 x \wedge B_2 x)$ where $\star$ is the Lorentzian Hodge star. 

\begin{lem}\label{lem:sigma}
%Suppose that $x\in \Om$. 
$\Si_{12}=\Si_{23}=\Si_{31}=\Si=C^+\setminus \Om$. 
%Then $x\in \Si_{12}$ implies that $x\in \Si_{13}$ and . 
\end{lem} 
\proof The following argument is borrowed from \cite{DU}. 
For $i\ne j$ consider the rational maps $\Phi_{ij}: [x]\in \P V \to \P (\La^2 V)$, where $V=\R^{n+1}$, given by projectivization of the correspondence
$$
x\mapsto B_i(x) \wedge B_j(x)
$$
It is clear that the domain of $\Phi_{ij}$ in $C^+$ is $C^+ \setminus \Si_{ij}$ and 
$\Phi_{ij}$ does not extend continuously to any point of $\Si_{ij}$. 
The assumption that the vectors $B_i(x), i=1,2,3$ are linearly dependent for all $x$ implies that $\Phi_{12}=\Phi_{23}=\Phi_{31}$. Therefore, their sets of indeterminacy $\Si_{12}, \Si_{23}, \Si_{31}$ are also equal. In particular, 
for every $x\in \Si$, $\rank(B_{\ul{A}}(x))=1$. \qed 

\medskip 
We now begin the actual proof of Theorem \ref{thm:simple}. 

\medskip
{\bf Case 1 (The generic case).} $A_i, i=1,2,3$ are all loxodromic and no two of them have a common eigenvector in $C$. 
In particular, in view of Lemma \ref{L3.1}, every eigenvector $x\in C$ of $A_i, i=1,2,3$, belongs to $\Om$. 

Let $x_1\in C$ be an eigenvector of $A_1$ with the eigenvalue $\la_1> 1$. Then for $x=x_1$, \eqref{eq:ld} implies: 
$$
[\al_1(x) (\la_1 -1) -\al_2(x) -\al_3(x)]x + \al_2(x) A_2(x) + \al_3(x) A_3(x)=0.  
$$
Note that all three vectors $x, A_2(x), A_3(x)$ belong to the cone $C$. Thus, the above equation implies that the vectors
$$
[\al_1(x) (\la_1 -1) -\al_2(x) -\al_3(x)]x ,\quad \al_2(x) A_2(x) ,\quad \al_3(x) A_3(x)
$$
belong to the same line. Since $x\in \Om$, the last two vectors are nonzero. 
%Note that we cannot have both $\al_2(x)=0, \al_3(x)=0$ since $\la_1>1$. 

\medskip
a. Suppose that $[\al_1(x) (\la_1 -1) -\al_2(x) -\al_3(x)]\ne 0$.  
Then $x$ is a common eigenvector for $A_2, A_3$ contradicting the assumptions of Case 1.  

b. Thus, $[\al_1(x) (\la_1 -1) -\al_2(x) -\al_3(x)]=0$ and
$$
\al_2(x) A_2(x) + \al_3(x) A_3(x)=0.
$$
Then, since $A_2(x), A_3(x)\in C$ (contained in the future light cone), it follows that 
$\al_2(x)$, $\al_3(x)$ have to have opposite signs.  
By applying the same argument to the eigenvectors $x_i$ of $A_i, i=2,3$, we obtain:
$$
\al_{i}(x_k)\al_j(x_k)<0,
$$
for all choices of pairwise distinct $i,j,k\in \{1,2,3\}$. 

It immediately follows that there is an index $i\in \{1,2,3\}$ such that the function $\al_i(x)$ changes its sign on the set $\Om\cap C$. 
However, $\Om$ is connected and $\al_i(x)\ne 0$ on $\Om$. Contradiction. 

\medskip
{\bf Case 2.} Suppose that $A_i, i=1,2,3$ are all loxodromic and $A_1, A_2$ share a common eigenvector in $C$. 
%, but $A_1, A_2, A_3$ does not generate a cyclic group. 

Then discreteness of $\Ga$ implies that $A_1, A_2$ share both eigenvectors in $C$. Let $P_{12}$ 
be the plane spanned by these eigenvectors.  If $\al_3(x)=0$ for some $x\in H\cap P_{12}$, 
then (by \eqref{eq:ld}) we get
$$
\al_1(x) A_1(x) +\al_2(x) A_2(x) = (\al_1(x) +\al_2(x)) x.
$$
By Lemma \ref{L2.0}, it follows that $A_1(x)=A_2(x)=x$, contradicting the assumption that $A_1, A_2$ are loxodromic. 
Thus, $\al_3(x)\ne 0$ for all $x\in H\cap P_{12}$. Then for all $x\in H\cap P_{12}$,  the equation \eqref{eq:ld} 
%$$\sum_{i=1}^3 \al_i(x) B_i(x)=0, $$
implies that $A_3(x)$ is a linear combination of $A_1(x), A_2(x), x$ which all belong to $P_{12}$. It then follows that $A_3$ preserves 
 $L=H\cap P_{12}$, i.e., $A_1, A_2, A_3$ have a common axis in the hyperbolic plane. 
 Consider the group $\<A_1, A_2, A_3\>$ generated by 
 $A_1, A_2, A_3$. This group acts discretely on $L$ (since is a subgroup of the discrete group $\Ga< PO(3,1)$). If the action of 
 $\<A_1, A_2, A_3\>$ on $L$ were not faithful, this group would contain an elliptic element fixing $L$ pointwise. This contradicts our 
 assumption that all elliptic elements of $\Ga$ are Cartan involutions. Hence, the group  $\<A_1, A_2, A_3\>$ 
 acts faithfully on $L$ as a discrete group of translations. Therefore, $\<A_1, A_2, A_3\>\cong \Z$, contradicting the hypothesis 
 of Theorem \ref{thm:simple}.

\medskip
{\bf Case 3.} $A_1$, $A_2$ are loxodromic and $A_3=J$ is elliptic (a Cartan involution).  
Let $p\in H$ be the unique fixed vector of $A_3$. Then 
$p=p_3\in \Si_{13}$ (see Lemma \ref{L3.1}). 
By Lemma \ref{lem:sigma}, it follows that $p\in \Si_{12}$. Since $p$ is not an eigenvector of $A_i, i=1,2,$ 
 it follows (by Lemma \ref{L3.1}) that $A_1(p)=A_2(p)$, i.e., $A_2^{-1} A_1(p)=p$. Thus, $A_2^{-1} A_1$ is elliptic fixing $p$. 
 Since $\Ga$ belongs to the class $\CK$, $A_2^{-1} A_1$ is a Cartan involution. 
 Therefore,  $A_2^{-1} A_1=J=A_3$ (since a Cartan involution is determined by its fixed point). 

Our goal is to obtain a contradiction. Let $x_{i}^{\pm}\in C$ be the eigenvectors of $A_i, i=1,2$ with eigenvalues $\la_i^{\pm 1}, i=1,2$. 
Let $x=x_{1}^{\pm}$. Note that $x$ cannot be an eigenvector of $J$.

\begin{lemma}\label{L3.2}
Either $x$ is an eigenvector of $A_2$, or $J(x)$ is a multiple of $A_2(x)$.  
%$i=1,2$ is taken modulo $2$. 
\end{lemma}
\proof Our proof is similar to the argument in Case 1. 
We will assume that $x$ is not an eigenvector of $A_2$. Then, by Lemma \ref{L3.1}, $x\in \Om$. 
In particular, $\al_i(x)\ne 0, i=1,2,3$. We have: 
$$
[\al_1(x) (\la_1^{\pm 1} -1) -\al_2(x) -\al_3(x)]x + \al_2(x) A_2(x) + \al_3(x) J(x)=0.  
$$
As before, the vectors $x, A_2(x), A_3(x)$ belong to the cone $C$. Thus, the vectors
$$
[\al_1(x) (\la_1 -1) -\al_2(x) -\al_3(x)]x ,\quad \al_2(x) A_2(x) ,\quad \al_3(x) A_3(x)
$$
belong to the same line and the last two vectors are nonzero. 
%we cannot have both $\al_2(x)=0, \al_3(x)=0$ since $\la_1\ne 1$. 

\medskip
a. Suppose that $[\al_1(x) (\la_1^{\pm 1} -1) -\al_2(x) -\al_3(x)]\ne 0$. 
Then $x$ is a common eigenvector of $A_2, A_3$ contradicting our assumption that $A_3=J$ is a Cartan involution.

b. Thus, $\al_2(x)\al_3(x)\ne 0$ and  $[\al_1(x) (\la_1^{\pm 1} -1) -\al_2(x) -\al_3(x)]=0$. 
Hence, $A_2(x), J(x)$ are multiples of each other. \qed  

The same argument, of course, applies to the eigenvectors of $A_2$. %and we obtain: 

\begin{cor}
 One of the following holds:

a. $A_1, A_2$ have a common axis and commute. 
%$A_{i+1}$ commutes with $A_i$ ($i$ is taken modulo $2$), or 

b. For each $i=1,2$, $C_i=A_{i+1}^{-1} J$ and $A_i$ generate a cyclic group ($i$ is taken modulo $2$). 
\end{cor} 
\proof If $A_1, A_2$ has a common eigenvector in $C$, then, by discreteness of $\Ga$, they share both eigenvectors in 
$C$ and, hence, have a common axis in $\H^3$. Since $A_1, A_2\in \Ga$ and $\Ga$ is in the class $\CK$, 
it follows that $A_1, A_2$ commute. Thus, suppose that $[A_1, A_2]\ne 1$. Then, by Lemma \ref{L3.2}, for 
$x=x_i^{\pm}$, $C_i=A_{i+1}^{-1}J(x)$ is a multiple of $x$. Hence, the elements $A_i, C_i$ share both eigenvectors 
in $C$. Therefore, they have the same axis in $\H^n$ and, since $\Ga$ is in the class $\CK$, these elements have to 
generate a cyclic group. \qed

\medskip 
Note that, since $A_3=J=A_2^{-1} A_1$, it follows that in the case (a) of this corollary, 
all three elements $A_1, A_2, A_3$ commute. This is impossible 
since $A_1, A_2$ are loxodromic and $J$ is a Cartan involution. Thus, (b) holds for both $i=1,2$ and  
$A_i, C_i=A_{i+1}^{-1} J$ generate a cyclic group. 

\medskip 
Combining the equations
$$
A_2= A_1J, \quad A_2= JC_1^{-1},
$$
we obtain
\begin{equation}\label{eq:con}
J A_1 J= C_1^{-1}. 
\end{equation}
Therefore, $J$ preserves the axis $L_1$ of $A_1$ in $\H^n$. Since $J$ is a Cartan involution, it has to reverse orientation on $L_1$. 
We write $A_1=A R$, where $A$ is a hyperbolic element with the axis $L_1$ and $R$ is an elliptic element fixing $L_1$ pointwise. In particular, $RJ=JR$. Then, $C_1^{-1}= A^{-1} R$ and $C_1=A R^{-1}=A_1 R^{-2}$ and
$$
A_1 C_1^{-1}= R^2. 
$$  
Since $A_1, C_1\in \Ga$, we also have $R^2\in \Ga$. By our assumptions on elliptic elements of $\Ga$, $R^2=1$. Thus, $C_1=A_1$. 
For the same reason, $C_2=A_2$. Hence, by equation \eqref{eq:con}, $J$ anticommutes with both $A_1, A_2$. 
In particular, the fixed point of $J$ belongs to  $L_1\cap L_2$ and $J$ preserves both $L_1$ and $L_2$. 

However, $A_2=A_1J$ and, since $A_1, J$ preserve $L_1$, it follows that $A_2$ also preserves $L_1$ as well, 
i.e., $A_1, A_2$ are loxodromic isometries with the common axis $L=L_1=L_2$. 
But then the composition  $J=A_2^{-1}A_1$ has to be either loxodromic or elliptic fixing $L$ or the identity. 
This contradicts the assumption that $J$ is a Cartan involution.

\medskip
{\bf Case 4.} $A_2, A_3$ are (distinct) elliptic of order $2$, so $A_i=A_i^{-1}, i=2,3$. 
(We make no assumptions about $A_1$ apart from $A_1\ne A_2, A_1\ne A_3$.) The same arguments as in Case 3 
(considering fixed points $p_2, p_3$ of $A_2, A_3$) show that $A_2 A_1= A_2^{-1} A_1=A_3, A_3 A_1= A_3^{-1} A_1= A_2$. Thus,
$$
A_2= A_1 A_3= A_3 A_1
$$
and, hence $A_1, A_3$ commute. Since $A_3$ is a Cartan involution and $A_1$ is loxodromic or Cartan, it follows that 
$A_1=A_3$, which contradicts the assumption that the elements $A_1, A_2, A_3$ are all distinct. 

This concludes the proof of Theorem \ref{thm:simple}. \qed

\section{Dirichlet domains of cyclic loxodromic groups}

In this section we discuss the exceptional case of triples of elements of cyclic loxodromic groups. 
The following result is implicit in the work of T.~Drumm and J.~Poritz \cite[\S5, \S 7]{Drumm-Poritz}, who analyzed Dirichlet fundamental domains of cyclic subgroups of $SO(3,1)^\uparrow$ in great detail: 

\begin{theorem}\label{thm:cyclic}
Let $\<A\> < SO(3,1)^\uparrow\cong PSL(2,\C)$ be a cyclic loxodromic group. Then the Dirichlet tiling $\CD_x$ of $\<A\>$ in $\H^3$ is simple for 
every choice of $x\in \H^3$. 
\end{theorem}

\begin{conjecture}
The same conclusion holds for all cyclic loxodromic subgroups of $PO(n,1), n\ge 3$. 
\end{conjecture}

{\em Proof of Theorem \ref{thm:cyclic}.} We will be using notation and terminology of \cite{Drumm-Poritz}. In particular, we will use the notation 
$X_n=Bis(x, A^n x), n\ne 0$, for the bisectors. We will be using the notation $\Delta:=D_x$ for a fixed Dirichlet domain, and 
$F_n$ for $X_n\cap \Delta$, provided that this intersection is 2-dimensional. We also use the notation $S_n$ for the intersections of the ideal boundaries of $X_n$ and $\Delta$, provided that this intersection is 1-dimensional (an edge of the circular polygon $\geo \Delta$). 

\medskip 
1. Let $v\in \H^3$ be a vertex of $\Delta=D_x$. According to the conclusion on the bottom of page 177 of \cite{Drumm-Poritz}, the 
vertex $v$ is {\em splendid}, i.e., it belongs to exactly three faces $F_i, F_j, F_{i+j}$ of $\Delta$. If $\CD_x$ is not simple at $v$, by Lemma \ref{L2.1} and the following remark, $v$ belongs to a bisector $X_n$ so that $X_n\cap \Delta$ is not a 2-face. By Proposition 7.6 of  \cite{Drumm-Poritz}, $X_n\cap \Delta\ne v$. By  Proposition 7.7 of \cite{Drumm-Poritz}, $X_n$ cannot contain a finite edge $E$ of $\Delta$ incident to $v$.  By 
Proposition 7.8 of \cite{Drumm-Poritz}, $X_n$ cannot contain an infinite edge $E$ of $\Delta$ incident to $v$. Thus, 
$\CD_x$ is simple at $v$. 
 
2. Let $E$ be a bi-infinite edge of $\Delta$ such that $\CD_x$ is not simple along $E$. By Lemma \ref{L2.1}, $E=X_n\cap \Delta$ for 
some bisector $X_n$. Let $v, w\in S^2$ denote the ideal boundary points of $E$: $\geo E=\{v,w\}$. 
Then the ideal boundary circle of $X_n$ passes through $v, w$. 
By Corollary 5.7 of  \cite{Drumm-Poritz}, it follows that the ideal boundary of $\Delta$ has exactly four sides and, by Lemma 5.5  of 
\cite{Drumm-Poritz} these sides are: $S_i, S_j, S_{-i}, S_{-j}$. Without loss of generality, we can assume that $v=S_i\cap S_j$. Then, by 
Corollary 5.7 of  \cite{Drumm-Poritz}, $n=i+j$. Up to relabeling, there are two options for the vertex $w$:

a. $w=S_{-i}\cap S_{-j}$. However, by Lemma 5.5  of \cite{Drumm-Poritz}, $w\in \geo X_{-i-j}\cap \geo \Delta$. Since the involution 
$\phi$ defined in \cite{Drumm-Poritz} swaps $\geo X_{i+j}\cap \geo \Delta$ and $\geo X_{-i-j}\cap \geo \Delta$, it follows that
$$
\geo X_{-i-j}\cap \geo \Delta= \{v,w\}. 
$$
By repeating the arguments in the proof of Proposition 4.4 of \cite{Drumm-Poritz}, we see that  $-i-j=n$. 
Hence, $-i-j=n=i+j$ and $n=0$, contradiction. 

b. $w=S_{-j}\cap S_{i}$. Then $w\in \geo X_{i-j}\cap \geo \Delta$. In this case there is no reason to expect that 
$\{v,w\}= \geo X_{i-j}\cap \geo \Delta$. Nevertheless, by Proposition 4.5 of \cite{Drumm-Poritz}, we get:
$$
A^{j-i}(w)\in A^{-n}(\{v,w\}). 
$$
If $A^{j-i}(w)= A^{-n}(w)$ then (as in the proof of  Proposition 4.4 of \cite{Drumm-Poritz}) $i+j=n=j-i$ and, 
hence, $i=0$, contradiction. If $A^{j-i}(w)= A^{-n}(v)=A^{-i-j}(v)$ then 
$$
A^{-2j}(v)=w. 
$$
However, by looking at the fundamental domain $\geo \Delta$ we also see that
$$
A^{-j}(v)=w,
$$
thus $-2j=-j$ and $j=0$. Contradiction. \qed

\medskip 
In view of Theorem \ref{thm:cyclic}, it remains to consider cyclic subgroups of $O(3,1)^\uparrow$ generated by 
orientation-reversing loxodromic isometries $A$. Such isometries $A$ are called {\em glide-reflections}:  
$A$ is the composition $A_0 R$, where $A_0$ is hyperbolic and $R$ is a reflection 
in a hyperplane containing the axis $L$ of $A_0$.  

\begin{proposition}\label{prop:cyclic}
For every $A$ as above, the Dirichlet tiling $\CD_x$ of $\<A\>$ in $\H^3$ is simple for 
every choice of $x\in \H^3$. 
\end{proposition}
\proof One can, in principle, go through the proofs given in \cite{Drumm-Poritz} and modify them when necessary in order to allow 
orientation-reversing loxodromic elements. Instead, we will give a direct argument. 
%We first observe that if $n\ne 0$ is even then $A^n$ is hyperbolic. 
%Furthermore, for distinct even $n, m$ the bisectors $Bis(x, A^n x), Bis(x, A^m x)$ are disjoint in $\H^3$. Therefore, we 
%need to consider triples $A_1, A_2, A_3\in \<A\>$ where at most one element (say, $A_2$) is hyperbolic. 

Let $L\subset \H^3$ denote the  axis of $A$. Since $A$ is orientation-reversing, there exists a hyperbolic plane $P\subset \H^3$ containing $L$, invariant under $A$, so that $A$ reverses orientation on $P$. Hence, $A$ preserves the half-spaces bounded by $P$. For a point 
$x\in \H^3$ let $x_P\in P$ denote the point nearest to $x$. The nearest-point projection $x\to x_P$ commutes with the action of $A$. 

\begin{lemma}\label{lem:perp bisectors}
For every $x\in \H^3$ and $n\in \Z\setminus \{0\}$, the bisector $Bis(x, A^n x)$ is orthogonal to $P$ and 
$$
Bis(x, A^n x)= Bis(x_P, A^n x_P). 
$$
\end{lemma}
\proof We set $y:=A^nx$ and let $p\in \R^{3,1}$ be such that $P=p^\perp \cap H$. 
The extended bisector $\widetilde{Bis}(x,y)\subset \R^{3,1}$ equals $(x-y)^\perp$. 
Computing $(x-A^n x)\cdot p$ and  using the fact that 
$Ap=p$, we obtain: $(x-A^n x)\cdot p= x\cdot p -x\cdot A^{-n} p=0$. Thus, $p\in \widetilde{Bis}(x,y)\subset \R^{3,1}$. Therefore, since 
$\widetilde{Bis}(x,y)=(x-y)^\perp$ and $P=p^\perp\cap H$, 
the hyperplanes $\be:=Bis(x,y)$ and $P$ in $\H^3$ are orthogonal. Let $R_\be\in PO(3,1)$ 
be the isometric reflection in the hyperplane $\be$. Since $\be$ is orthogonal to $P$, $R_\be$ preserves $P$. In particular, 
$R_\be$ commutes with the projection $z\to z_P, z\in \H^3$. Since $R_\be(x)=y$, it follows that $R_\be(x_P)=y_P$. Therefore, 
$\be=Bis(x,y)$ 
is the bisector for $x_P, y_P=A^n x_P$ as well. \qed 

\medskip 
In view of this lemma, Dirichlet tilings $\CD_x$ and $\CD_{x_P}$ (with respect to the group $\<A\>$) are the same. 
Therefore, it suffices to prove simplicity of the Dirichlet tilings $\CD_{x_P}$ 
of $\<A\>$ on $P=\H^2$. We, thus, assume that $x\in P$. The isometry $A$ of $\H^2$ is the composition of the hyperbolic isometry $A_0$ preserving $L$ and the reflection $R$ in $\H^2$ fixing $L$.

Let $x_L\in L$ denote the point in $L$ nearest to $x$. Again, the nearest-point projection $x\to x_L$ commutes with the action of $A$. For $n\in \Z\setminus \{0\}$ we let $m_n\in L$ denote 
the midpoint of $x_L, A^n x_L$. We claim that $Bis(x, A^n x)$ passes through $m_n$. Indeed, similarly to Lemma \ref{lem:perp bisectors},
$$
Bis(A_0^n x, x)= Bis(A^n x_L, x_L).
$$
Hence, $d(x, m_n)=d(A_0^nx, m_n)$, while 
$$d(m_n, A_0^nx)= d(m_n, R^n A_0^n x)= d(m_n, A^n x).$$
Thus, $d(x, m_n)=d(A^nx, m_n)$ and $m_n\in Bis(x, A^nx)$ proving the claim. 

We now consider the bisectors $Bis(x, A^{\pm 1} x), Bis(x, A^{\pm 2} x)$. These bisectors bound a convex polygon $F\subset \H^2$ (of infinite area) containing $x$. The vertices of this polygon are $y= Bis(A^{-2}x, x)\cap  Bis(A^{-1}x, x)$, 
$z=Bis(A^{-1}x, x)\cap  Bis(Ax, x)$ and $w= Bis(A^{2}x, x)\cap  Bis(Ax, x)$. We next observe that
$$
d(y,L)=d(z,L)=d(w,L). 
$$
This follows from congruence of the triangles 
$$
\Delta(y m_{-2} m_{-1}), \quad \Delta(m_{-1} x_L z), \quad  \Delta(m_{1} x_L z), \quad \Delta(w m_{2} m_{1}).
$$ 
Since $A$ sends $Bis(A^{-1}x, x)$ to  $Bis(Ax, x)$ and preserves the distance to $L$, it follows that
$$
y \stackrel{A}{\longrightarrow} z \stackrel{A}{\longrightarrow} w
$$
and $A^2: y\to w$. Thus, the polygon $F$ is a fundamental domain for the action of the group $\<A\>$ on $\H^2$. 
Since, clearly, $D_x\subset F$, we have $D_x=F$. Furthermore, the only vertex-cycle of the fundamental domain $F$ is 
$$
y \stackrel{A}{\longrightarrow} z \stackrel{A}{\longrightarrow} w \stackrel{A^{-2}}{\longrightarrow} y,
$$
which has length $3$. Therefore, in the Dirichlet tiling $\CD_x$ in $\H^2$, there are exactly three fundamental tiles adjacent to each of the $y, z, w$. Hence, $\CD_x$ is simple. Proposition \ref{prop:cyclic} follows. \qed 

\begin{figure}[htbp] %  figure placement: here, top, bottom, or page
   \centering
   \includegraphics[width=3in]{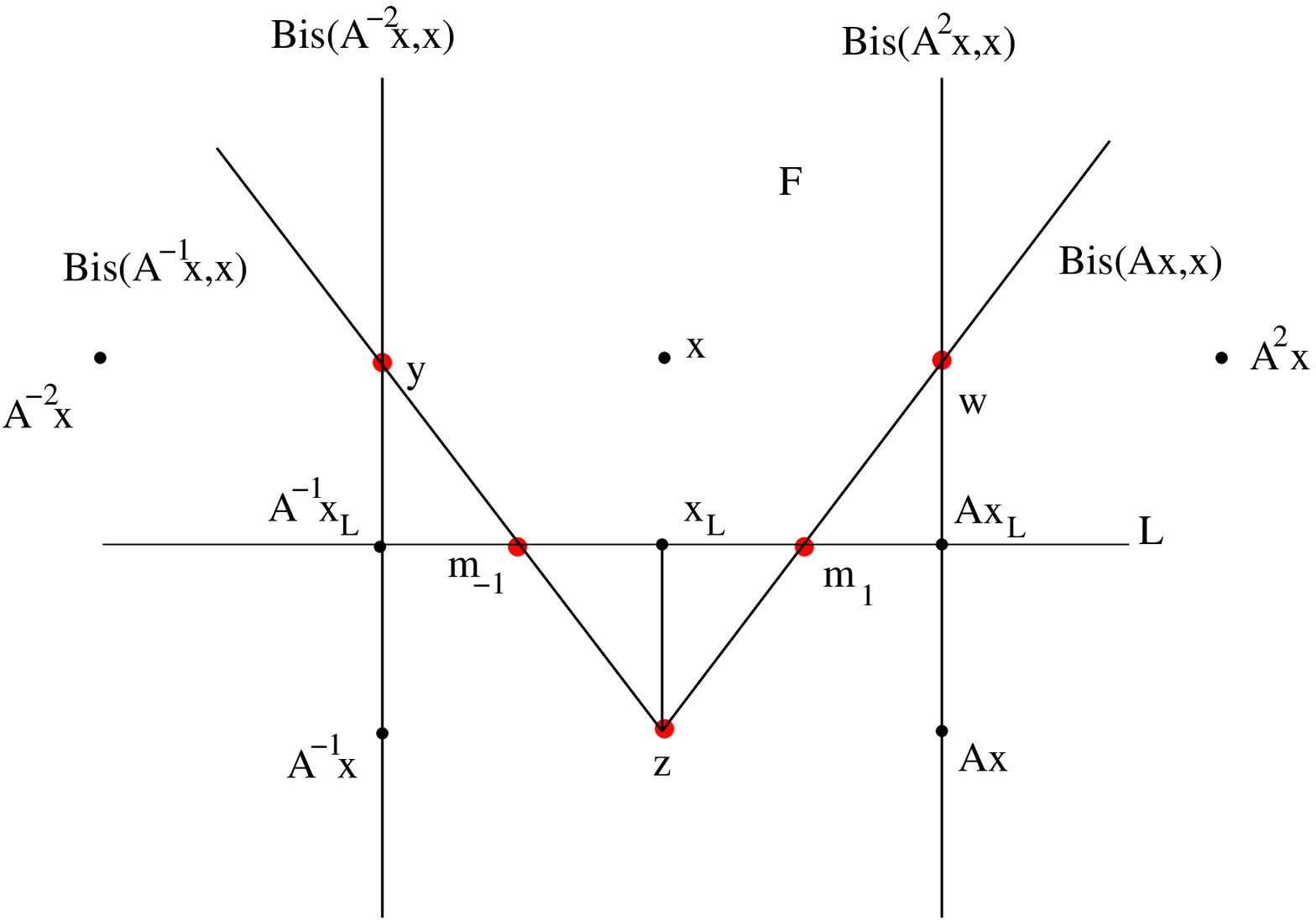} 
      \label{dir.fig}
\end{figure}

\section{Two examples}\label{sec:examples}

\begin{example}
There exists a cyclic loxodromic subgroup $\<A\> < SO(3,1)^\uparrow$ for which there exists an open nonempty 
subset $U\subset \H^3$ so that for all $x\in U$ the triple intersection of bisectors 
$$
Bis(A^{-1}x, x)\cap Bis(A^2 x, x)\cap Bis(A^3 x, x)\subset \H^3
$$
is non-transversal (i.e., is a hyperbolic geodesic).  
\end{example}

Let $A\in PO(2,1)$ be an orientation-reversing loxodromic isometry of $\H^2$. We extend $A$ to an orientation-preserving isometry 
of $\H^3$ (also denoted $A$). We will consider a triple of distinct 
nontrivial elements $A_1, A_2, A_3\in \<A\>$ such that $A_1, A_3$ are orientation-reversing and $A_2$ is orientation-preserving isometries of $\H^2$. 

Let $B: \R^4\to Mat_{4,3}$ be the associated mapping $x\mapsto (B_1(x), B_2(x), B_3(x))$, where $B_i=A_i-I$, see equation \eqref{eq:B}. 

\begin{lemma}
$B$ is singular if and only if $A_2=A_1 A_3$. 
\end{lemma}
\proof We choose the basis of eigenvectors $e_1, e_2, e_3, e_4$ of $A$ in $\R^4$, where $e_1, e_2\in C$ are eigenvectors with 
eigenvalues $\la, \la^{-1}\ne 1$, $e_1\cdot e_2=-1$, the unit vectors $e_3, e_4$ are orthogonal to $e_1, e_2$ and each other and
$$
Ae_i=-e_i, i=3,4. 
$$
We now consider vectors $x=(x_1,x_2,x_3,x_4)\in \R^4$ so that $x_1 x_2 x_3x_4\ne 0$. Let $\la_i, \la_i^{-1}$ be the eigenvalues of 
$A_i$ corresponding to the eigenvectors $e_1, e_2$ respectively. Consider the $4\times 3$ matrix $B(x)$, $x\in \R^4$. 
The two bottom rows of this matrix are $x_i[-2, 0, -2]$, $i=3, 4$. Hence, $\rank(B(x))= 2$ if and only if the following determinant equals zero: 
$$
\Delta= 
\left| 
\begin{array}{ccc}
\la_1 -1 & \la_2 -1       & \la_3 -1\\
\la_1^{-1} -1 &\la_2^{-1} -1 & \la_3^{-1} -1\\
1                 & 0                 & 1
\end{array}
\right|.$$
(Hence, this is independent of $x$.) Computing $\Delta$ we obtain: 
$$ 
\Delta=\left|
\begin{array}{cc}
\la_1 -1 & \la_2 -1 \\
\la_1^{-1} -1 &\la_2^{-1} -1 
\end{array}
\right| + 
\left| 
\begin{array}{cc}
\la_2 -1       & \la_3 -1\\
\la_2^{-1} -1 & \la_3^{-1} -1
\end{array}
\right|= 
\left| 
\begin{array}{cc}
\la_2 -1       & \la_3- \la_1\\
\la_2^{-1} -1 & \la_3^{-1}- \la_1^{-1}
\end{array}
\right|. 
$$
This determinant equals $0$ iff
$$
\frac{\la_2 -1}{\la_2^{-1} -1}=  \frac{\la_3 -\la_1}{\la_3^{-1} - \la_1^{-1}} \iff
$$
$$
\la_2 = \la_1 \la_3. 
$$
Equivalently, $A_2=A_1 A_3$. \qed

\begin{remark}
By adopting arguments from \cite{Drumm-Poritz}, one can prove that the same conclusion holds for all triples of 
loxodromic elements $\ul{A}=\{A_1, A_2, A_3\}$ of $PO(n,1)$ generating a cyclic group: 
After reordering these elements if necessary,  $A_2=A_1 A_3$ if and only if the associated map $B_{\ul{A}}$ is singular.  
\end{remark}

\medskip
To get a specific example, we will take $A_1=A^{-1}$, $A_2=A^2, A_3=A^{3}$. Our next goal is to find conditions on $x$ and $\la$ under which the Gramm matrix $Gr(x)$ of the vectors $\{B_1(x), B_2(x)\}$ is positive-definite, i.e., when the restriction of the Lorentzian inner product to 
$$
Span(B_1(x), B_2(x))^\perp
$$
is indefinite, that is, $Bis(A_1x, x)\cap Bis(A_2 x, x)\cap H\ne \emptyset$. We have: 
$$
Gr(x)= 
\left[
\begin{array}{cc}
-2(\la-1)(\la^{-1}-1)x_1 x_2 +4 x_3^2 					&   -\mu  x_1 x_2            \\
-\mu    x_1 x_2	 &  -2(\la^2-1)(\la^{-2}-1) x_1 x_2
  \end{array}
  \right] 
$$ 
where
$$
\mu= (\la-1)(\la^2-1)  + (\la^{-1}-1)(\la^{-2}-1) =  (\la-1)(\la^2-1) (1+ \la^{-3}).      
$$
Then 
$$
\det(Gr(x))= - \nu^2 x_1^2 x_2^2 + 8 (\la^{2}-1)^2 \la^{-2} x_1 x_2 x_3^2,
$$
where
$$
\nu = (\la-1)(\la^2-1) - (\la^{-1}-1)(\la^{-2} -1) = (\la-1)(\la^2-1) (1- \la^{-3})
$$
In addition, we have the condition $x\in \CL^\uparrow$, i.e., $x_3^2< 2x_1 x_2$. We now fix $x\in H$ 
such that $x_3\ne 0$ and let $\la\to 1+$. Then, 
$$
\left( (\la-1)(\la^2-1) (1- \la^{-3}) \right)^2 \sim (\la-1)^6 = o( (\la^{2}-1)^2 \la^{-2} )
$$
as $\la\to 1+$. This means that each all $\la>1$ sufficiently close to $1$, $det(Gr(x))>0$ and, hence, 
the open set 
$$
U_\la=\{x\in H: \det(Gr(x))>0\}
$$
is nonempty. Hence, 
$Bis(A_1x, x)\cap Bis(A_2 x, x) \cap Bis(A_3 x, x)\cap H\ne \emptyset$ for all $x\in U_\la$. 
The reader who enjoys computations will verify that for every $\la$ with
$$
1< \la < \frac{3+\sqrt{5}}{2}\approx 2.6,
$$
the set $U_\la$ is nonempty. Therefore, for all such $\la$, there exists an open nonempty set $U_\la\subset \H^3$ so that for 
all $x\in U_\la$ 
$$
Bis(A_1x, x)\cap Bis(A_2 x, x)\cap Bis(A_3 x, x)\subset \H^3$$
is a complete hyperbolic geodesic, i.e., the triple intersection of bisectors in $\H^3$ is nontransversal. Furthermore, 
the loxodromic elements $A_i$ belong to a cyclic group $\<A\>$ of orientation-preserving isometries of $\H^3$ that stabilize 
a hyperbolic plane $\H^2\subset \H^3$. (The group $\<A\>$ does not preserve the orientation on $\H^2$.) \qed 

\medskip
Note that the example constructed above does not contradict Theorem \ref{thm:cyclic}: The nontransversal 
triple intersections do not occur on the boundary of the Dirichlet domain $D_x$. In our second example of a discrete abelian subgroup of $SO(3,1)$, such nontransversal intersections occur on the boundary of $D_x$ for an open nonempty set of $x\in \H^3$. 

\begin{example}
Let $A$ be a hyperbolic isometry of $\H^3$ with the axis $L$  
and let $R\in PO(3,1)$ be the order 2 elliptic rotation around $L$.  We consider the abelian 
group $\Ga:=\<A, R\>$ generated by these isometries. Let $A_1:=R, A_2:=A, A_3:=RA$. 
Then for all $x\in \H^3\setminus L$ the triple intersection of bisectors 
$$
Bis(A_1x, x)\cap Bis(A_2 x, x)\cap Bis(A_3 x, x)\subset \H^3
$$
is nonempty. Furthermore, this line of intersection $I_x$ is contained in the boundary of the Dirichlet domain $D_x$ of $\Ga$. 
\end{example}

Note that in this example, $\Ga$ preserves hyperbolic planes $P\subset \H^3$ containing $L$ and $x$ ($\Ga$ reverses the orientation 
on $P$). We first compute the fundamental domain $D_x$ for $\Ga$: It is the solid $S$ in $\H^3$ containing $x$ and 
bounded by the bisectors $Bis(Ax, x), Bis(A^{-1} x), Bis(x, Rx)$. Indeed, clearly, $S$ contains $D_x$. On the other hand, $S$ is a fundamental polyhedron for $\Ga$ which can be easily verified using the Poincar\'e's fundamental domain theorem, see e.g. \cite{Maskit}, \cite{Ratcliffe}. Since $D_x$ is also a fundamental domain of $\Ga$, it follows that $S=D_x$. Next, the intersection  
$$
I_x:= Bis(Ax, x)\cap Bis(R x, x)\subset \H^3
$$
is a hyperbolic geodesic in $\H^3$ contained in the boundary of $D_x$ (since the bisectors $Bis(Ax, x), Bis(A^{-1} x)$ are disjoint). 
Set $Q:=Bis(R x, x)$. Then the reflection $R_Q$ in the hyperplane $Q$ sends $x$ to $Rx$ and $Ax$ to $RAx$ and fixes $I_x$. 
Therefore, for every $y\in I_x$, $d(y,x)=d(y, Ax)=d(y, RAx)$. Hence, 
$$
I_x= Bis(Ax, x)\cap Bis(R x, x)\cap Bis(RAx, x). \qed 
$$

\section{Proof of Theorem \ref{thm:generic}}\label{sec:generic}

We now can prove Theorem \ref{thm:generic}. In view of Lemma \ref{L2.1}, we need to prove that for generic choice of $x\in \H^3$, 
for every edge $E\subset D_x=D_x(\Ga)$, $E$ is the intersection of exactly two bisectors $Bis(A_1x, x), Bis(A_2,x, x)$, where 
$A_1, A_2\in \Ga$. First, for every triple $\ul{A}=\{A_1, A_2, A_3\}$ of nontrivial distinct elements of $\Ga$ which do not belong 
to a common cyclic subgroup, we define the set $\CE(\ul{A})$ consisting of those $x\in \R^4$ for which the intersection
$$
\bigcap_{i=1}^3 \widetilde{Bis}(A_i x, x)\subset \R^4
$$
has dimension $\ge 2$. In other words, 
$$
\CE(\ul{A})=\{x\in  \R^4: \rank( B_{\ul{A}}(x) )\le 2 \}, 
$$
see \S \ref{sec:complexes} for the notation. 

The set $\CE(\ul{A})$ is clearly algebraic in $\R^4$ and is stable under multiplication of $x$ by scalars.  
According to  Theorem \ref{thm:simple}, $\CE(\ul{A})$ is a proper algebraic subset of $\R^4$. 
In particular, its intersection with $\CL^\uparrow$ is closed and has topological dimension $\le 3$. Since $\CE(\ul{A})$ is stable under 
scaling, the intersection $\CE_H(\ul{A}):=\CE(\ul{A})\cap H$ is nowhere dense. 
Since $\Ga$ is countable, the union $\CE_H$ of the subsets $\CE_H(\ul{A})$ (taken over all 
triples $\ul{A}$ of distinct nontrivial elements of $\Ga$ generating non-cyclic groups)  
is nowhere dense in $H=\H^3$. Next, consider the triples 
$\{A_1, A_2, A_3\}$ of distinct  nontrivial elements of $\Ga$ generating cyclic subgroups. For such a triple, 
by   Theorem \ref{thm:cyclic} and Proposition \ref{prop:cyclic}, 
$$
\bigcap_{i=1}^3 Bis(A_i x, x) \subset \H^3
$$
is disjoint from $D_x(\<A\>)$ for every choice of $x\in \H^3$. Since $D_x(\Ga)\subset D_x(\<A\>)$, it follows that such 
nontransversal triple intersection is disjoint from $D_x(\Ga)$ as well, so we can ignore such triples. Thus, we conclude that $\CD_x$ 
is weakly simple for all $x\in \H^3\setminus  \CE_H$.

It remains now to show that for generic $x$, fixed points of 
Cartan involutions in $\Ga$ do not belong to any edge of $\CD_x$.  

\begin{lemma}
Let $p\in H$, $A\in O(3,1)^\uparrow$ is an element not fixing $p$. Then there exists a hyperbolic plane $Q_p\subset H$ such that for 
all $x\in H\setminus Q_p$, $p\notin Bis(Ax, x)$.  
\end{lemma}
\proof $p\in Bis(Ax, x)$ if and only if: 
$$
p\cdot (Ax- x) =0 \iff (Ap -p)\cdot x=0$$ 
 Since $Ap\ne p$, the orthogonal complement to the vector $Ap-p$ is a hyperplane in $\R^4$. For every $x$ away from this hyperplane, 
 $p\notin Bis(Ax, x)$. \qed 

\medskip
Since $\Ga$ is countable, $\H^3$ contains only countably many fixed points $p_i, i\in \N,$ 
of Cartan involutions that belong to $\Ga$. Therefore, for 
every 
$$
x\in \H^3 \setminus (\CE_H \cup \bigcup_{i\in \N} Q_{p_i}),
$$
for every Cartan involution $J\in \Ga$, the fixed point $p$ of $J$ in $H$ belongs to the unique bisector $Bis(Ax, x), A\in \Ga$, 
namely, $Bis(Jx,x)$. Hence, $p$ cannot belong to an edge of $\CD_x(\Ga)$. Theorem \ref{thm:generic} follows. \qed 

\section{3-dimensional hyperbolic orbifolds}

The goal of this section is to prove Theorem \ref{T1.1}. Our proof is a minor variation of the one in \cite{PP}.

We recall that 
an {\em orbihedron} $O$ is a topological space $|O|$ (the {\em underlying space} of $O$) together with an atlas where each chart is the quotient $U_\al/G_\al$ of a polyhedral complex $U_\al$ by a finite PL group action $G_\al\acts U_\al$, 
satisfying certain compatibility conditions, see e.g. \cite{Haefliger}. 
An orbihedron is called an {\em orbifold} if the polyhedral complexes above are PL manifolds. The {\em singular set} $\Si_O$ of an orbihedron is the subset of $|O|$ consisting of points $x$ which are covered by $\tilde{x}\in U_\al$ with nontrivial stabilizer in 
$G_\al$. The {\em order} of a singular point $x$ is the order of the stabilizer of $\tilde{x}$ in $G_\al$. 
An orbifold is called a DISK, an ANNULUS or a TORUS, if it is the quotient of a disk or an annulus or a torus by a finite group action. 
(See \cite{Kapovich2000}.) For instance, the Moebius band is an ANNULUS.

\begin{notation}
Suppose that $S$ is a surface. We let $S(m_1,...,m_k)$, where $m_i=2,3,...,\infty$ denote the 2-dimensional orbifold with 
boundary obtained from $S$ as follows: 

1. For each $i$ with $m_i=\infty$, we remove the interior of a closed disk from $S$, so that the disks are pairwise disjoint. 

2. For each $i$ so that $m_i<\infty$, we introduce a singular point of order $m_i$ on $S$ (away from the disks removed and from the boundary of $S$). 

In order to shorten the notation, if $m_1=m_2=...=m_\ell=m$, we replace the repeating sequence $(m_1,...,m_\ell)$ in our notation with 
$\ell \times m$.  For instance, $$S^2(\infty, \infty, 2,2,2,2)=S^2(2\times \infty, 4\times 2),$$ 
the annulus with four singular points of order $2$.  
\end{notation}

Below we will use the notation $I$ for the interval $[-1,1]$. 

\medskip 
We now review the construction given in \cite{PP}. 
Define a regular 2-dimensional cell complex $X$ obtained from $\R \P^2$ by attaching 2-cells 
$D_1, D_2$ to two distinct projective lines $L_1, L_2$ in $\R \P^2$. The lines $L_1, L_2$ cut $\R\P^2$ in two 2-cells $D_3, D_4$. Next, as in \cite{PP}, define the 2-dimensional orbihedron $Y$ by introducing 3 singular points of order $2$ in the interior of 
each of the disks $D_i, i=1,...,4$. It is proven in \cite{PP} that for every finitely-presented group $G$ there exists a finite orbi--cover $\tilde{Y}\to Y$, such that $\pi_1(|\tilde{Y}|)\cong G$.

Panov and Petrunin in \cite{PP} then ``thicken'' $Y$ to a compact 3-dimensional hyperbolic orbifold $Y_3$ with convex boundary, where each singular point $p_j$ of $Y$ corresponds to a singular segment $p_j\times I$ (of the order $2$) in $Y_3$; in addition, $Y_3$ constructed in \cite{PP} has an extra order $2$ singular point $q_i$ for each thickened disk $D_i, i=1,2$. Then, $Y_3$ constructed in \cite{PP} is the quotient of a closed convex   subset $C\subset \H^3$ by a discrete convex-cocompact group of isometries $\Ga_3<PO(3,1)$, which contains both orientation-preserving elliptic involutions (corresponding to the 
singular segments $p_j\times I$) and Cartan involutions (corresponding to the isolated singular points $q_i)$. 

We now observe that instead of thickening the orbihedron $Y$ described above, we can thicken a slightly different one: Let $Y^+$ be the 
orbihedron obtained from $Y$ by adding an extra order 2 singular point $q_i$ to each cell $D_i, i=3,4$. 
(Nothing changes as far as the 2-cells $D_1, D_2$ are concerned, they still have three order 2 singular points each.) 
Now,  if $f: \tilde{Y}\to Y$ is an orbi-cover, it induces an orbi-cover 
$f^+: \tilde{Y}^+\to Y^+$, which is unramified over the points $q_3, q_4$: The orbihedron $\tilde{Y}^+$ is obtained from $\tilde{Y}$ by 
declaring the points in $f^{-1}(q_3)\cup f^{-1}(q_4)$ to be singular points of order $2$. Clearly, 
$\pi_1(|\tilde{Y}^+|)=\pi_1(|\tilde{Y}|)$. Therefore, as in \cite{PP}, for every finitely-presented group $G$ there exists a finite orbi-cover $\tilde{Y}^+\to Y^+$, such that $\pi_1(|\tilde{Y}^+|)\cong G$. We will see below why thickening the orbihedron 
$Y^+$ is better than thickening $Y$. 

Before thickening $Y^+$, we describe its double cover (as in \cite{PP}), 
which will, hopefully, clarify the construction. Let $\Pi: S^2\to \R \P^2$ be the 2-fold cover, quotient by the antipodal involution 
$\tau$. We let $\al_i:=\t{L}_i, i=1,2, \t{D}_j, j=3,4$ denote the complete preimages of the lines and the disks under $\Pi$. Note that each $\t{D}_3$ and $\t{D}_4$ consists of two copies of $D_3$ and $D_4$ respectively, while $\al_1, \al_2$ are circles. 
We lift the orbifold data accordingly, so we obtain the orbifold $O^2$ (sphere with 16 singular points of order 2) 
which is the 2-fold cover of $\R \P^2(8\times 2)$. 
Now, thickening $\al_1\cup \al_2$ in $S^2$ yields a quadruply-punctured sphere $F$. %(This is the same 4-holed sphere which appears below.) 
Of course, the restriction $\Pi: \al_i\to L_i$ is again a 2-fold cover. The 2-fold cover $S^1\to S^1$ extends to a 2-fold orbifold-cover $T^2(\infty)\to D^2(3\times 2)$ (see below). Thus, attaching   two copies of $T^2(\infty)$ (the one-holed torus) 
to $O^2$ along the loops $\al_1, \al_2$, we obtain an orbihedron $\hat{Y}^+$ which is a 2-fold orbi-cover of $Y^+$. 

Instead 
of thickening $Y^+$ we will equivariantly thicken its 2-fold cover $\hat{Y}^+$: The surface $F$ is thickened to $F\times I$, while 
both copies of $T^2(\infty)$ are thickened to $T^2(\infty)\times I$. The 3-manifolds $Z_i:=T^2(\infty)\times I$ are then attached to 
$F\times I$ along the appropriate annuli ${A}_{\al_i},i=1,2$, in $F\times \{\pm 1\}$ 
(thickenings of the loops $\al_1\times \{1\}$ and $\al_2\times \{2\}$), which are identified with the annuli 
$\D T^2(\infty)\times I$. Lastly, the four orbi-disks  in 
$$
\t{D}_3\cup \t{D}_4
$$
will be thickened to the appropriate orbifold $I$-bundles ${W}_i, i=1,2,3,4$, over $D^2(4\times 2)$ and attached to 
$F\times I$ along $\D F\times I$. (The precise construction of ${W}_i$ will be given below.) 
It is then clear (e.g., from Van Kampen's theorem) that the fundamental group of the 
resulting orbifold $\hat{O}$ is isomorphic to $\pi_1(\hat{Y}^+)$. Assuming that $\tau$ extends to an involution of 
$\hat{O}$ (with isolated fixed-points only), we obtain the 3-dimensional orbifold $\CO=\hat{O}/\tau$. 
By the construction, $\CO$ and $Y^+$ have isomorphic fundamental groups. 

\medskip 
We now explain how to construct ${W}_i$'s and how to extend the involution $\tau$.  
Begin with the 2-torus $T^2$ and its elliptic involution $\si: T^2\to T^2$: 
It has 4 fixed points and the quotient $T^2/\si$, as an orbifold, is $S^2(4\times 2)$. We extend $\si$ to the orientation--reversing involution
$$
\si: T^2\times I\to T^2\times I, \quad \si(z,t)=(\si(z), -t). 
$$
The orbifold $(T^2\times I)/\si$ has only isolated singular points (four of them). 
Then the projection $\eta: (T^2\times I)/\si\to T^2/\si$ is the orbifold $I$-bundle. This projection is the quotient of the projection 
if $T^2\times I$ to the 1st factor. 

\medskip
{\bf Definition of $W_i$'s.} 
We define $W$ to be the suborbifold of $(T^2\times I)/\si$ obtained by removing $\eta^{-1}(D)\cong D\times I$, 
where $D$ is a   nonsingular 2-disk in $S^2(4\times 2)$. In particular,  
$\eta^{-1}(\D D)$ is an annulus. Then the orbifolds $W_i, i=1,...,4$ are copies of $W$ above. They will be attached to $F\times I$ 
by gluing the annuli $\eta^{-1}(\D D)$ to the annuli in $\D F\times I$. 

\medskip
{\bf Extension of $\tau$.} We extend $\tau$ to $F\times I$ by the identity to the second factor. Then $\tau$ sends each 
$A_{\al_i}$ to itself, where $A_{\al_i}\subset F \times \D I$ is an annular thickening of $\al_i\times \{\pm 1\}$. The quotient 
$A_{\al_i}/\tau$ is the Moebius band (since $\tau$ reverses orientation on $F$). We are identifying   
$A_{\al_i}$ with the annulus $A_i\subset \D Z_i$, where $Z_i= T^2(\infty) \times I$ and $A_i=\D T(\infty)\times I$. 
Note that $\tau$ acts on the annulus $A=A_i=S^1 \times I$ by $\tau(z,t)= (\tau(z), -t)$, where $\tau: S^1\to S^1$ is an involution. 
(We now drop the index $i$ since the construction is the same for $i=1$ and $i=2$.)  
Thus, we again take the elliptic involution $\si: T^2\to T^2$. Let $x\in T^2$ be one of its four fixed points. Take a small 
$\si$-invariant   2-disk $D\subset T^2$ around $x$. We then regard $T^2(\infty)$ as $T^2\setminus int(D)$. The involution $\si$ restricts to the involution $S^1\to S^1$ of the boundary circle of   $T^2(\infty)$ which is isotopic to $\tau: S^1\to S^1$, so we identify them. Set $Z:=T^2(\infty)\times I$. Now, the map $\si: Z\to Z$ given by 
$$
\si(z,t)=(\si(z), -t)
$$
is the required extension of $\tau$ to $Z=Z_i, i=1,2$. Clearly, the orbifold 
$$
V=Z/\tau$$
 has only 3 singular points. This concludes the construction of $\hat{O}$ and the extension $\tau: \hat{O}\to \hat{O}$. Therefore, we obtain the 3-dimensional orbifold with boundary 
$\CO:= \hat{O}/\tau$ which is a thickening of the orbihedron $Y^+$. Furthermore, the singular locus of $\CO$ is finite. We also have the orbifold-fibration 
$$\zeta: V\to D^2(2,2,2)$$ obtained as the quotient (by $\tau$) of the projection $Z\to T^2(\infty)$ to the first factor.

\medskip
{\bf Topological properties of $\CO$}. Our goal is to show that the orbifold $\CO$ is {\em hyperbolizable}, i.e., 
there exists a closed  convex $Q$ of $\H^3$ and a discrete isometry group $\Ga_{\CO}< PO(3,1)$, so that the quotient-orbifold $Q/\Ga_{\CO}$ is homeomorphic to the orbifold $\CO$. In principle, this could be proven by constructing $\Ga_{\CO}$ by hand, via Maskit combination. 
Instead, we will show that $\CO$ is hyperbolizable by verifying that it is irreducible and atoroidal, in which case $\CO$ is 
hyperbolizable by Thurston's hyperbolization theorem, see e.g., \cite{Bonahon, Kapovich2000}. 

\medskip 
We first analyze the JSJ decomposition of the orbifold $\CO$, see e.g. \cite{Bonahon}. 
Recall  that $\CO$ is constructed 
from 5 pieces: Orbifold $N:=(F\times I)/\tau$, two copies of the orbifold  $V:=Z/\tau$ (where $Z=T^2(\infty)\times I$) and two copies of the orbifold 
$W= (T^2(\infty,\infty) \times I)/\tau$. We now convert each of these orbifolds to an orbifold pair by marking 
some of their boundary annuli/Moebius bands:

1. Define $(N, P_N)$, where $N=F\times I$ and $P_N=  \partial F \times I \sqcup A_{\al_1}\times \{1\} \sqcup A_{\al_2}\times \{-1\}$. Then 
set $(U, P_U):=(N/\tau, P_N/\tau)$. 

2. Define $(V, P_V)$, where $V:=Z/\tau$ and $P_V= \zeta^{-1}(\D D^2(2,2,2))$ is a single Moebius band. 

3. Define $(W, P_W)$, where $P_W=  \eta^{-1}(\D D)$ is a single annulus, see above for the definition of $\eta: (T^2\times I)/\si\to T^2/\si$.

\medskip
For each of the orbifolds $U, V, W$ we define its {\em partial boundary} $\D_P$ by: $\D_P U:= \D U\setminus P_U$, etc.

By the construction, each of the orbifolds $N, W, V$ is very good: It admits a finite manifold-cover. Also, each of the orbifolds is strongly atoroidal, i.e., it contains no $\pi_1$-injective TORI: 
Its fundamental group is virtually free and, hence, contains no $\Z^2$. 
%Note also that the orbifold-pairs $(V, P_V)$ and 
%$(W, P_W)$ admit a common finite cover since  this is true for 2-dimensional hyperbolic orbifolds with boundary  
%$D^2(2,2,2)$ and $D^2(2,2,2,2)$. Since being atoroidal and irreducible is preserved by finite covers, we will consider only the pair 
%$(V, P_V)$. 

\begin{lemma}
The orbifold pairs $(U\setminus \D_P U, P_U)$, $(V\setminus \D_V ,P_V)$ and $(W\setminus \D_W ,P_W)$ are all irreducible, boundary-irreducible and 
acylindrical, see \cite{Kapovich2000} for the terminology). 
\end{lemma}
\proof We will give a proof for $(V\setminus \D_V, P_V)$ since the other pairs are similar. We first note that irreducibility and acylindricity are stable under passing to finite covers. Now, $V$ is finitely covered by the product $M:=T^2(\infty)\times I$, so 
that $P_V$ lifts to the annulus $P_M:=\D T^2(\infty) \times I$. Irreducibility of $M$ is clear. Boundary-irreducibility follows from the fact that the annulus $P_M$ is $\pi_1$-injective in $M$. To see that $(M\setminus \D_P M, P_M)$ is acylindrical, note that the 
image of $\pi_1(P_M)$ in $\pi_1(M)\cong \Z \star \Z$ is a maximal cyclic subgroup of $\pi_1(M)$. \qed

\begin{lemma}
The orbifold-pair $(U, P_U)$ is acylindrical.  
\end{lemma}
\proof It suffices to prove acylindricity for the 2-fold cover $(N, P_N)$ of $(U, P_U)$. Then surface $P_N$ contains the union of 
the annuli in $\D F\times I$. Every annulus properly embedded in $F\times I$ and  disjoint from  
$\D F\times I$, is isotopic to one of the form $c\times I$, where $c$ is a simple loop in $F$. On the other hand, every essential simple 
loop $c$ in $F$ (i.e., a loop which does not bound a disk and does not bound an annulus whose other boundary component is in $\D F$) 
has to cross either $\al_1$ or $\al_2$. Therefore, the corresponding annulus $c\times I$ either crosses $A_{\al_1}$ or $A_{\al_2}$. Hence, 
it cannot be isotopic to the annulus disjoint from the union of circles $\D P_N$. Thus, $(U, P_U)$ is acylindrical. \qed 

\medskip
We now conclude that the sub-orbifolds $W_i, i=1,2$ and $V_i, i=1,2$ are maximal (up to isotopy) $I$-bundles in $\CO$. Therefore, their union is the {\em characteristic suborbifold} in $\CO$, and, hence, its splitting along the annuli and Moebius bands $P_{V_i}, P_{W_j}$ 
is its JSJ decomposition. Irreducibility of $\CO$ follows from the fact that each annulus and Moebius band 
$P_{V_i}, P_{W_j}$ is incompressible in $\CO$ (as it is incompressible in the pieces of the JSJ decomposition). In particular, $\CO$ contains no bad suborbifolds. Atoroidality of $\CO$ follows since each essential TORUS in $\CO$ has to be contained in one of the characteristic 
suborbifolds and they are all strongly atoroidal. We thus proved:

\begin{proposition}
The orbifold $\CO$ is hyperbolizable: It can be realized as the orbifold-quotient of a closed convex subset $Q$ of $\H^3$ by a discrete 
isometry group $\Ga_{\CO}<PO(3,1)$. In particular, the group $\Ga_{\CO}$ is a convex-cocompact subgroup of $PO(3,1)$. 
\end{proposition}

We then observe that $\Ga_{\CO}\cong \pi_1(\CO)$ contains a free nonabelian subgroup, say, $\pi_1(F)$. In particular, the group 
$\Ga_{\CO}$ is nonelementary. We can now finish the proof of Theorem \ref{T1.1}. Given a finitely-presented group $G$ we find 
a finite index subgroup $\t\Ga< \pi_1(Y^+)=\Ga_{\CO}$ so that 
$$
G\cong \t\Ga/\<\<torsion\>\>. 
$$
 The group $\t\Ga$ is the fundamental group of some orbifold $\hat\CO$ (a finite covering of $\CO$). Since $\CO$ is hyperbolizable, 
 we obtain a discrete embedding $\t\Ga \embed PO(3,1)$. Since all singularities of $\CO$ are isolated, so are all singularities of 
 its finite cover $\hat\CO$. Thus, $\t\Ga$ belongs to class $\CK$. Since $\Ga_{\CO}$ is convex-cocompact and $\t\Ga$ has finite index in 
 $\Ga_{\CO}$, it follows that  $\t\Ga$  is also convex-cocompact. Thus, $\t\Ga$ belongs to the class $\CK^2$. Theorem \ref{T1.1} follows. \qed

%\begin{proposition}
%For every finitely-presented group $G$ there exists a finite orbifold 
%cover $\tilde{O}\to O$, such that $\pi_1(|\tilde{O}|)\cong G$.
%\end{proposition}

\section{Constructing projective varieties}

{\em Proof of Theorem \ref{thm:main}.} 
Let $G$ be a finitely-presented group. By Theorem \ref{T1.1}, there exists a nonelementary group $\t\Ga<PO(3,1)$ of class $\CK^2$, so that 
$G\cong \pi_1(\H^3/\t\Ga)$. We let $x\in \H^3$ be a generic base-point, so that the associated 
Dirichlet tiling $\CD_x(\t\Ga)$ of $\H^3$ is weakly simple. Recall that since $\t\Ga$ is convex-cocompact, every face of $\CD_x(\t\Ga)$ 
is a finitely-sided convex hyperbolic polytope. However, since $\H^3/\t\Ga$ has infinite volume, so is the fundamental domain $D_x$. Therefore, unlike in \cite{KK}, $D_x$ is not a polytope from the projective viewpoint: It is bounded by linear subspaces as well as the quadric $\P C$. 
We will see, nevertheless, that this is harmless. 

%In principle, one can remedy this by replacing 
%$D_x$ (and polytopes in its $\Ga$-orbit) with their ``polytopal hulls'' in $\R\P^3$, where we 

We now define the locally finite hyperbolic polyhedral complex $\tilde\CC=\CD_x(\t\Ga) - \CD_x^{(0)}(\t\Ga)$, 
the derived complex of $\CD_x(\t\Ga)$ (see \S \ref{sec:complexes}). 
Let $\Ga< \t\Ga$ be a torsion--free normal finite index subgroup in $\t\Ga$, so that $\CC:=\tilde{\CC}/\Ga$ is a simple finite hyperbolic 
polyhedral complex, see Lemma \ref{lem:regular}.  Set $\Theta:=\t\Gamma/\Gamma$. This finite group acts naturally on $\CC$ and this action is transitive 
on facets (since the action of $\t\Ga$ is transitive on facets $D_{\ga x}$ of $\CD_x(\t\Ga)$). 

Consider the manifold $M:=\H^3/\Gamma$, and let $F\subset M$ denote the finite set which is the image of $\CD_x^{(0)}(\t\Ga)$ in $M$. 
If $m\in F$ is a vertex of $\CD_x(\t\Ga)/\Ga$, it is not fixed by any nontrivial element of $\Theta$ (since it is so for vertices of the complex 
$\CD_x(\t\Ga)$). Therefore, %for $E=E'/\Phi$, we have:
$$
\pi_1((M \setminus F)/\Theta)\cong \pi_1(M/\Theta)= \pi_1(\H^3/\t\Ga), 
$$
which is the quotient of $\t\Ga$ by the normal closure of the Cartan involutions in $\t\Ga$.

We next {\em complexify} the tiling $\CC$ as in \S \ref{sec:varieties}. The result is a complex $\CP$ of projective spaces based on $\CC$. The action 
$\Theta\acts \CC$ lifts to the action $\Theta\acts \CP$. 

\begin{lemma}
The action $\Theta\acts \CP$ satisfies the Assumption \ref{asm:cartan} in \S \ref{sec:varieties}. 
\end{lemma}
\proof As in \ref{sec:varieties}, we consider the complex of projective spaces $\tilde\CP$ based on the complex 
$\t{\CC}$. Each stratum $\tilde{X}_i$ is the span of the corresponding face $\tilde{c}_i$ of $\t{\CC}$ and 
$\tilde{X}_i$ projects isomorphically to the corresponding stratum $X_i$ of $\CP$. Therefore, it suffices to verify that the action 
$\t\Ga\acts \tilde\CP$ satisfies the Assumption \ref{asm:cartan}. Suppose that $\ga\in \t\Ga$ fixes a point $p\in \tilde{X}$, the push-out of $\tilde{\CP}$. 
Let $\tilde{X}_i$ be the minimal stratum of $\tilde{X}$ containing $p$. Then $\ga$ preserves $\tilde{X}_i$ and, hence, preserves the corresponding 
face $\tilde{c}_i$ of $\tilde\CC$. Thus, $\ga$ is elliptic and has to be a Cartan involution since $\t\Ga$ is in the class $\CK$. Hence, the fixed-point 
set of $\ga$ in  $\tilde{X}_i=\Span(\tilde{c}_i)$ is the disjoint union of the point $p$ and the dual (with respect to the Lorentzian inner product) 
projective space $p^\perp\subset \Span(\tilde{c}_i)$. Since $p^\perp$ is disjoint from $\H^3$, it is also disjoint from  
$\tilde{c}_i$. \qed

We next replace $\CP$ with its blowup $b\CP$ and let $X$ denote the projective variety which is the push-out of $b\CP$, see  \S \ref{sec:varieties}. 
All singularities of $X$ are normal crossings. The group 
$\Theta$ acts naturally on $X$ and this action is transitive on top-dimensional strata (since $\Theta$ acts transitively on facets of $\CC$). 
In particular, $Z:=X/\Theta$ is irreducible. 

Note that since $\Ga$ is nonelementary and torsion-free it  has trivial  center, see Lemma \ref{lem:nocenter}. 
Thus, by Proposition \ref{prop:centerless}, 
$$
\pi_1(Z)\cong \pi_1(C/\Theta)\cong \pi_1((M\setminus F)/\Theta)\cong \pi_1(\H^3/\t\Ga)=G. 
$$
All singularities of $Z$ are normal crossings and $\Z_2$-quotients of normal crossing singularities, types $Y_1$ and $Y_2$ described in \S \ref{sec:varieties}. 
By blowing up centers of type 1 singularities $Y_1$ and dividing by $\Theta$, we get a new irreducible projective variety $V$ where all singularities 
are normal crossings and their quotients of the type $Y_2/\Z_2$. As before, $\pi_1(V)\cong   G$. Lastly, by the argument in the end of \S \ref{sec:varieties}, 
we replace the 3-dimensional $V$ with irreducible projective surface $W$ so that $\pi_1(W)\cong \pi_1(V)\cong   G$ and all 
all singularities of $W$ are normal crossings and Whitney umbrellas.

In order to prove the second assertion of Theorem \ref{thm:main}, we note that if $\t\Ga$ is a torsion-free convex-cocompact subgroup 
of $PO(3,1)$, then the group $\Theta=\t\Ga/\Ga$ acts freely on $X$ and, hence, all singularities of $Z=X/\Theta$ are normal crossings. 
Then one takes $V:=Z$ and proceeds as above. This concludes the proof of Theorem \ref{thm:main}. \qed

\begin{remark}
It was proven by Carlson and Toledo in \cite{CT} that if $G$ is a K\"ahler group (e.g., fundamental group of a smooth projective variety) 
which is isomorphic to a nonelementary discrete subgroup $\Ga$ of $PO(n,1)$, then $G$ contains a finite index subgroup isomorphic to the fundamental group of a Riemann surface.  (Note that \cite{CT} assumes that $\Ga$ is cocompact, but it is clear from the proof that nonelementary is enough.) 
\end{remark}

\medskip
\noindent Address: Michael Kapovich: Department of Mathematics, University of California, Davis, CA 95616,
USA. (kapovich@math.ucdavis.edu)

\end{document}